\documentstyle[12pt,fullpage,emlines,bezier]{article}

\linethickness{0.4pt}

\newcommand{\bi}{\bibitem}
\newcommand{\nb}{\newblock}

\newcommand{\be}[1]{\begin{equation}\label{#1}}
\newcommand{\ee}{\end{equation}}

\newcommand{\la}{\langle\,}
\newcommand{\ra}{\,\rangle}
\newcommand{\ve}{\varepsilon}
\newcommand{\prf}{{\bf Proof.}\ }
\newcommand{\aaa}{{\cal A}}
\newcommand{\ccc}{{\cal C}}
\newcommand{\pp}{{\cal P}}
\newcommand{\rr}{{\cal R}}
\newcommand{\dd}{{\cal D}}
\newcommand{\bdelta}{\bar\delta}
\newcommand{\distan}{\mathop{\rm dist}}
\newcommand{\ppi}{{\cal P}=\la x\mid x^2=x\ra}
\newcommand{\topp}{\mathop{\mbox{\bf top}}}
\newcommand{\bott}{\mathop{\mbox{\bf bot}}}

\newtheorem{thm}{\quad Theorem}
\newtheorem{lm}{\quad Lemma}

\title{On the Properties of the Cayley Graph of Richard Thompson's Group $F$}
\author{\vspace{2ex}
V.~S.~Guba\thanks{This research is partially supported by the RFFI
grant 99--01--00894 and the INTAS grant 99--1224.}\\
Vologda State Pedagogical University,\\
6 S.\,Orlov Street,\\
Vologda\\
Russia\\
160600\\
E-mail: guba{@}uni-vologda.ac.ru}

\begin{document}

\maketitle

\begin{abstract}

We study some properties of the Cayley graph of the R.\,Thompson's group $F$
in generators $x_0$, $x_1$. We show that the density of this graph, that is,
the least upper bound of the average vertex degree of its finite subgraphs
is at least $3$. It is known that a $2$-generated group is not amenable if
and only if the density of the corresponding Cayley graph is strictly less
than $4$. It is well known this is also equivalent to the existence of a
doubling function on the Cayley graph. This means there exists a mapping
from the set of vertices into itself such that for some constant $K>0$,
each vertex moves into the distance at most $K$ and each vertex has at
least two preimages. We show that the density of the Cayley graph of a
$2$-generated graph does not exceed $3$ if and only if the group satisfies
the same condition with $K=1$. Besides, we give a very easy formula to find
the length (norm) of a given element of $F$ in generators $x_0$, $x_1$. This
simplifies the algorithm by Fordham. The length formula may be useful to
find the general growth function of $F$ in generators $x_0$, $x_1$ and the
growth rate of this function. In this paper we show that the lower bound for
the growth rate of $F$ is $(3+\sqrt5)/2$.

\end{abstract}

\section*{Introduction}

The Richard Thompson group $F$ can be defined by the following infinite
group presentation

\be{xinf}
\la x_0,x_1,x_2,\ldots\mid x_j{x_i}=x_ix_{j+1}\ (i<j)\,\ra.
\ee
This group was found by Richard~J.~Thompson in the 60s. We refer to the
survey~\cite{CFP} for details. (See also~\cite{BS,Bro,BG}.) It is easy
to see that for any $n\ge2$, one has $x_n=x_0^{-(n-1)}x_1x_0^{n-1}$ so
the group is generated by $x_0$, $x_1$. It can be given by the following
presentation with two defining relations

\be{x0-1}
\la x_0,x_1\mid x_1^{x_0^2}=x_1^{x_0x_1},x_1^{x_0^3}=x_1^{x_0^2x_1}\ra,
\ee
where $a^b=b^{-1}ab$ by definition.

Each element of $F$ can be uniquely represented by a {\em normal
form\/}, that is, an expression of the form
\be{nf}
x_{i_1}x_{i_2}\cdots x_{i_s}x_{j_t}^{-1}\cdots x_{j_2}^{-1}x_{j_1}^{-1},
\ee
where $s,t\ge0$, $0\le i_1\le i_2\le\cdots\le i_s$, $0\le j_1\le j_2
\le\cdots\le j_t$ and the following is true: if~(\ref{nf}) contains
both $x_i$ and $x_i^{-1}$ for some $i\ge0$, then it also contains
$x_{i+1}$ or $x_{i+1}^{-1}$ (in particular, $i_s\ne j_t$).

Equivalent definition of $F$ can be done in the following way. Let us
consider all strictly increasing continuous piecewise-linear functions
from the closed unit interval onto itself. Take only those of them that
are differentiable except at finitely many dyadic rational numbers and
such that all slopes (derivatives) are integer powers of $2$. These
functions form a group under composition. This group is isomorphic to $F$.
Another useful representation of $F$ by piecewise-linear functions can be
obtained if we replace $[0,1]$ by $[0,\infty)$ in the previous definition
and impose the restriction that all functions on infinity have the form
$t\mapsto t+c$, where $c$ is an integer.

The group $F$ has no free subgroups of rank $>1$. It is known that $F$ is
not elementary amenable (EA). However, the famous problem about amenability
of $F$ is still open. If $F$ is amenable, then it is an example of a
finitely presented amenable group, which is not EA. If it is not
amenable, then this gives an example of a finitely presented group,
which is not amenable and has no free subgroups of rank $>1$.  Note that
the first example of a non-amenable group without free non-abelian
subgroups has been constructed by Ol'shanskii~\cite{Olsh}. (The
question about such groups was formulated in~\cite{Day}, it is also
often attributed to von Neumann~\cite{vNeu}.) Adian~\cite{Ad83}
proved that free Burnside groups with $m>1$ generators of odd exponent
$n\ge665$ are not amenable. The first example of a finitely presented
non-amenable group without free non-abelian subgroups has been recently
constructed by Ol'shanskii and Sapir~\cite{OlSa}. Grigorchuk~\cite{Gri}
constructed the first example of a finitely presented amenable group
not in EA.
\vspace{2ex}

The author thanks Matt Brin and Goulnara Arjantseva for helpful
remarks.

\section{Density}
\label{dens}

By the {\em density\/} of a finite graph $\Gamma$ we mean the average
value of the degree of a vertex in $\Gamma$. More precisely, let $v_1$,
\dots, $v_k$ be all of vertices of $\Gamma$. Let $\deg_\Gamma(v)$ denote
the degree of a vertex $v$ in the graph $\Gamma$, that is, the number of
oriented edges of $\Gamma$ that come out of $v$. Then

\be{dgrform}
\delta(\Gamma)=\frac{\deg_\Gamma(v_1)+\cdots+\deg_\Gamma(v_k)}k
\ee
is the density of $\Gamma$.

Let $G$ be a group generated by a finite set $A$. Let $C(G,A)$ be the
corresponding (right) Cayley graph. Recall that the set of vertices of
this graph is $G$ and the set of edges is $G\times A^{\pm1}$. For an
edge $e=(g,a)$, its initial vertex is $g$, its terminal vertex is $ga$,
and the inverse edge is $e^{-1}=(ga,a^{-1})$. The {\em label\/} of $e$
equals $a$ by definition. For the Cayley graph $C=C(G,A)$ we define
the number

\be{denscayley}
\bdelta(C)=\sup\limits_\Gamma\delta(\Gamma),
\ee
where $\Gamma$ runs over all finite subgraphs of $C=C(G,A)$. So this number
is the least upper bound of densities of all finite subgraphs of $C$.
If $C$ is finite, then it is obvious that $\delta(C)=\bdelta(C)$. So we
may call $\bdelta(C)$ the {\em density of the Cayley graph\/} $C$.
\vspace{1ex}

Recall that a group $G$ is called {\em amenable\/} whenever there exists
a finitely additive normalized invariant mean on $G$, that is, a mapping
$\mu\colon{\cal P}(G)\to[0,1]$ such that $\mu(A\cup B)=\mu(A)+\mu(B)$
for any disjoint subsets $A,B\subseteq G$, $\mu(G)=1$, and $\mu(Ag)=
\mu(gA)=\mu(A)$ for any $A\subseteq G$, $g\in G$. One gets an equivalent
definition of amenability if only one-sided invariance of the mean is
assumed, say, the condition $\mu(Ag)=\mu(A)$ ($A\subseteq G$, $g\in G)$.
The proof can be found in~\cite{GrL}.

The class of amenable groups includes all finite and abelian groups. It is
invariant under taking subgroups, quotient groups, group extensions, and
ascending unions of groups. The closure of the class of finite and abelian
groups under these operations is the class EA of {\em elementary amenable\/}
groups. A free group of rank $>1$ is not amenable. There are many useful
criteria for (non)amenability~\cite{Fol,Kest,Gri80}. We need to mention the
two properties of a finitely generated group $G$ that are equivalent to
non-amenability.
\vspace{1ex}

{\bf NA$_1$.}\ {\sl If $G$ is generated by $m$ elements and $C$ is the
corresponding Cayley graph, then the density of $C$ has no maximum
value, that is, $\bdelta(C)<2m$.}
\vspace{1ex}

Note that if NA$_1$ holds for at least one finite generating set,
then the group is not amenable and so the same property holds for
any finite generating set. For the proof of this property, we need
to use the well-known {\em F\o{}lner condition\/}~\cite{Fol}. For
our reasons it is convenient to formulate this condition as follows.

Let $C$ be the Cayley graph of a group. By $\distan(u,v)$ we denote the
distance between two vertices in $C$, that is, the length of a shortest
path in $C$ that connects vertices $u$, $v$. For any vertex $v$ and a
number $r$ let $B_r(v)$ denote the ball of radius $r$ around $v$, that is,
the set of all vertices in $C$ on the distance at most $r$ from $v$. For
any set $Y$ of vertices, by $B_r(Y)$ we denote the $r$-neighbourhood of $Y$,
that is, the union of all balls $B_r(v)$, where $v$ runs over $Y$. By
$\partial Y$ we denote the {\em boundary\/} of $Y$, that is, the set
$B_1(Y)\setminus Y$. The F\o{}lner condition (for the case of a finitely
generated group) says that $G$ is amenable whenever $\inf\#\partial Y/\#Y=0$,
where the infimum is taken over all non-empty finite subsets in $G$ in
a Cayley graph of $G$ in finite number of generators (this property does not
depend on the choice of a finite generating set). Any finite set $Y$ of
vertices in $C$ defines a finite subgraph (also denoted by $Y$). The degree
of any vertex $v$ in $C$ equals $2m$, where $m$ is the number of generators.
We know that exactly $\deg_Y(v)$ of the $2m$ edges that come out of $v$,
connect the vertex $v$ to a vertex from $Y$. The other $2m-\deg_Y(v)$ edges
connect $v$ to a vertex from $\partial Y$. Note that each vertex of
$\partial Y$ is connected by an edge to at least one vertex in $Y$. This
implies that the cardinality of $\partial Y$ does not exceed the sum
$\sum(2m-\deg_Y(v))$ over all vertices of $Y$. Dividing by $\#Y$ (the number
of vertices in $Y$) implies the inequality $\#\partial Y/\#Y\le2m-\delta(Y)$.
If $\bdelta(C)=2m$, then $Y$ can be chosen such that $\delta(Y)$ is
arbitrarily close to $2m$ so $\#\partial Y/\#Y$ will be arbitrarily close to
$0$. On the other hand, for any vertex $v$ in $Y$ there are at most $2m$
edges that connect $v$ to a vertex in $Y$. Therefore, the sum
$\sum(2m-d_Y(v))$ does not exceed $2m\#\partial Y$. So
$2m-\delta(Y)\le2m\#\partial Y/\#Y$. If the right hand side can be done
arbitrarily close to $0$, then $\delta(Y)$ approaches $2m$ so $\bdelta(C)=2m$.
\vspace{1ex}

{\bf NA$_2$.}\ {\sl If $C$ is the Cayley graph of $G$ in a finite set of
generators, then there exists a function $\phi\colon G\to G$ such that
a$)$ for all $g\in G$ the distance $\distan(g,\phi(g))$ is bounded by a
constant $K>0$, b$)$ any element $g\in G$ has at least two preimages
under $\phi$.}
\vspace{1ex}

An elegant proof of this criterion based on the Hall-Rado theorem can be
found in~\cite{CGH}, see also~\cite{DeSS}. Note that this property also does
not depend on the choice of a finite generated set. A function $\phi$ from
NA$_2$ will be called a {\em doubling function\/} on the Cayley graph $C$.

We need a definition. Suppose that NA$_2$ holds for the Cayley graph
of a group $G$ for the case $K=1$. Then we say that the Cayley graph
$C$ is {\em strongly non-amenable\/}. The function $\phi\colon G\to G$ will
be called a {\em strong doubling function\/} on the Cayley graph $C$.
Note that each vertex is either invariant under $\phi$ or it maps into
a neighbour vertex. We know that NA$_2$ holds if and only if the group
is not amenable, that is, $\bdelta(C)<2m$. Now we would like to find
out what happens if the Cayley graph of a $2$-generated group is strongly
non-amenable.

\begin{thm}
\label{denle3}
The Cayley graph of a group with two generators is strongly non-amenable
if and only if the density of this graph does not exceed $3$.
\end{thm}

\prf Let $G$ be a group with $2$ generators and let $C$ be the corresponding
Cayley graph. It follows from the proof of~\cite[Theorem 32]{CGH} that $C$
is strongly non-amenable if and only if the doubling inequality holds,
that is, $\#B_1(Y)\ge2\#Y$ for any finite set $Y$ of vertices in $C$. Indeed,
if $C$ admits a doubling function $\phi$ from NA$_2$ with $K=1$, then
$\phi^{-1}(Y)$ has at least $2\#Y$ elements and it is contained in $B_1(Y)$.
To prove the converse, one needs to consider a bipartite graph with two
classes of vertices both equal to $G$. Two vertices from the different
classes are connected by an edge whenever the distance between them in $C$
does not exceed $K=1$. The doubling inequality implies that the conditions
of the Hall-Rado theorem hold. Therefore, the bipartite graph has a
perfect $(2,1)$-matching. This means that $C$ admits a doubling function
with $K=1$.

The ``only if" part is trivial since the doubling inequality is equivalent
to the fact that $\#\partial Y\ge\#Y$ for any finite subset $Y$ in $C$.
We know that $\#\partial Y/\#Y\le2m-\delta(Y)$ for any non-empty finite
subset $Y\subseteq C$, where $m=2$. This means that $\delta(Y)\le2m-1=3$
and so $\bdelta(C)\le3$.

Now suppose that $\bdelta(C)\le3$. Let $Y$ be a finite subgraph of $C$ with
$k$ vertices. For each $0\le s\le4$, let $q_s$ be the number of vertices in
$Y$ that have degree $s$ in $Y$. Clearly, $k=q_0+q_1+q_2+q_3+q_4$. The
number of oriented edges in $Y$ is the sum of all degrees of vertices, that
is, $q_1+2q_2+3q_3+4q_4$. If we divide this number by $k$, then we get
$\delta(Y)$. Since $\bdelta(C)\le3$, we have $\delta(Y)\le3$, which is
equivalent to the inequality $q_1+2q_2+3q_3+4q_4\le3k=3(q_0+q_1+q_2+q_3+q_4)$
and so it can be rewritten as $3q_0+2q_1+q_2-q_4\ge0$. For any finite
subgraph $Y$, we denote the number $3q_0+2q_1+q_2-q_4$ by $q(Y)$.

Suppose that the doubling inequality does not hold. Let us choose a minimal
counterexample $Y$ to the doubling inequality, that is, a finite subgraph $Y$
in $C$ with the property $\#B_1(Y)<2\#Y$ for which the number $q(Y)\ge0$ takes
the smallest possible value. Each vertex in $\partial Y=B_1(Y)\setminus Y$ is
connected by an edge with at least one vertex in $Y$. Suppose that some
vertex $v$ in $\partial Y$ is connected with at least two vertices in $Y$.
Let $Y'$ be the subgraph in $C$ with the new vertex $v$ added to $Y$ and two
new non-oriented edges that connect $v$ with vertices $v_1$, $v_2$ in $Y$.
By definition, $B_1(Y')=B_1(Y)\cup B_1(v)$. The vertex $v$ has exactly $4$
edges that come out of $v$. At least two of them connect $v$ with a vertex
in $Y$. So $B_1(v)$ may contain at most $2$ vertices not in $B_1(Y)$ (the
vertex $v$ itself belongs to $B_1(Y)$). This means that
$\#B_1(Y')\le\#B_1(Y)+2<2\#Y+2=2\#Y'$. Hence $Y'$ is also a counterexample.
To complete the proof, we need to check that $q(Y')<q(Y)$.

When we add a non-oriented edge that connects $v$ and $v_j$ ($j=1,2$), then
the degree of $v_j$ increases by $1$. This means that if $\deg_Y(v_j)=i$,
then $\deg_{Y'}(v_j)=i+1$. So $0\le i<4$, the value of $q_i$ decreases
by $1$, the value of $q_{i+1}$ increases by $1$. Clearly,
$q(Y)=3q_0+2q_1+q_2-q_4$ decreases by $1$. So if we add the edges for
both $v_1$ and $v_2$ (this does not exclude the case $v_1=v_2$), then
$q(Y)$ decreases by $2$. But we also have a new vertex $v$ that has degree
$2$ in $Y'$. Thus $q_2$ increases by $1$ and so after all these operations
we have $q(Y')=q(Y)-1$. This contradicts the minimality of $Y$.

The proof is complete.
\vspace{1ex}

Note that the proof of Theorem~\ref{denle3} goes without any changes
if we apply it to any regular graph of degree $4$. Also we have to
mention that the density of a Cayley graph of a group is closely related
to an isoperimetric constant $\iota_*$ of a graph (see the definition
in~\cite{CGH}). Namely, one has the equality $\iota_*(C)+\bdelta(C)=2m$
for the Cayley graph $C$ of an $m$-generated group.

Theorem~\ref{denle3} applied to the Cayley graph $\ccc_2$ of $F$ in
generators $x_0$, $x_1$ means that if we cannot find a subgraph in
$\ccc_2$ with density greater than $3$, then there exists a doubling
function on $\ccc_2$. One can imagine this doubling function in the
following way. Suppose that a bug lives in each vertex of $\ccc_2$. We
allow these bugs to jump at the same time such that each bug either
returns to its initial position or it jumps to a neighbour vertex. As
a result, we must have at least two bugs in each vertex.

It is natural to ask how much the value of $\delta(Y)$ can be for the
finite subgraphs we are able to construct. It is easy to see that in
each finite subgraph $Y$ of $\ccc_2$ there exists a vertex of degree at
most $2$. Indeed, if all vertices were of degree $3$ or $4$, then one could
travel along $Y$ by positively labelled edges only (if we enter a vertex
by an edge labelled by $x_0$ or $x_1$, then we can leave this vertex
travelling along an edge with one of these labels). But $F$ does not
have nontrivial relations that involve positive letters only.

However, the following result shows that we are able to construct
finite subgraphs in $\ccc_2$ with the density arbitrarily close to $3$.

\begin{thm}
\label{close3}
Let $\ccc_2$ be the Cayley graph of $F$ in generators $x_0$, $x_1$.
For any integers $m,n\ge1$ there exists a finite subgraph $\bar\Gamma_{n,m}$
in $\ccc_2$ such that
\be{limit}
\lim\limits_{m\to\infty}\delta(\bar\Gamma_{n,m})=\frac{6(n-1)}{2n-1}
\ee
for any $n\ge2$. In particular, $\bdelta(\ccc_2)\ge3$.
\end{thm}

\prf It is convenient to work with subgraphs in $\ccc=\ccc_\infty$. Note
that all Cayley graphs $\ccc_i$ ($i=2,3,\dots$) are subgraphs in $\ccc$.
Let $Y$ be a finite connected subgraph in $\ccc$. Assume that $Y$ is
{\em full\/}, that is, if two vertices of $Y$ can be connected in $\ccc$ by
an edge, then this edge belongs to $Y$. It is not hard to see that any
finite subset in $F$ determines the corresponding full subgraph, which
will be finite.

Let $Y$ be a finite connected full subgraph in $\ccc$. By a {\em rank\/} of
a vertex $v$ in $Y$ we mean the maximum number $k\ge0$ such that $v$ has an
edge in $Y$ labelled by $x_k^{\pm1}$ that comes out of $v$. Suppose that for
some $i\ge0$, the rank of any vertex in $Y$ is greater than $i$. Under these
conditions, we define a subgraph in $\ccc$ denoted by $\aaa_i(Y)$.

Let $v_1$, \dots, $v_k$ be the list of all vertices in $Y$ (they are
elements of $F$). Let $r_1$, \dots, $r_k$ be the ranks of these vertices
(each of these numbers exceeds $i$). For any $1\le j\le k$, we consider the
set of vertices of the form $v_jx_i^{-s}$, where $0\le s<r_j-i$. We call
this set a {\em column\/} of $v_j$. The number of vertices in this column
equals $r_j-i$ (if $r-j=i+1$, then the column consists of $v_j$ only).
The vertices in a column are connected in $\ccc$ by edges with label
$x_i$. The number of these edges in the column of $v_j$ equals $r_j-i-1$.
Note that all the columns are disjoint. The union of them will be the set
of vertices of the graph $\aaa_i(Y)$ we want to define.

Now let us take an arbitrary edge from $Y$. Suppose that it connects
$v'$ and $v''$. Let $r'$, $r''$ be the ranks of these vertices, respectively.
Without loss of generality we may assume that our edge has label $x_k$,
where $k>i$ so we have $v'x_k=v''$. If $k=i+1$, then we do nothing with this
edge. If $k=i+2$, then $v'x_i^{-1}$ and $v''x_i^{-1}$ both appear in the
columns of $v'$ and $v''$. In the graph $\ccc$, these vertices are
connected by an edge labelled by $x_{i+1}$ since $v'x_i^{-1}x_{i+1}=
v'x_{i+2}x_i^{-1}=v''x_i^{-1}$. In general, if $k=i+p+1$, where $p\ge1$,
then the rank of each of $v'$, $v''$ is at least $i+p+1$. So the elements
$v'x_i^{-t}$, $v''x_i^{-t}$ will appear in the columns for all $1\le t\le p$.
From the defining relations of $F$ it follows that $v'x_i^{-t}$,
$v''x_i^{-t}$ are connected in $\ccc$ by an edge with label $x_{k-t}$.
Indeed, $v'x_i^{-t}x_{k-t}=v'x_kx_i^{-t}=v''x_i^{-t}$ (we used the fact
$k-t>i$ and the relations of the form $x_i^{-1}x_j=x_{j+1}x_i^{-1}$
that hold in $F$ for all $j>i$).

So for any edge $e$ labelled by $x_k$ from $v'$ to $v''$ in $Y$ ($k\ge i+1$),
we have $k-i$ edges (including $e$) that connect a vertex in the column of
$v'$ with a vertex in the column of $v''$. These edges are called
{\em parallel\/} to $e$. The union of all these edges over all edges $e$
in $Y$ forms the set of edges of $\aaa_i(Y)$.

So we have defined the new graph $\aaa_i(Y)$. It is finite, connected,
and it contains $Y$. Let us check that $\aaa_i(Y)$ is full. Suppose that
two vertices $v'x_i^{-s}$ and $v''x_i^{-t}$ of $\aaa_i(Y)$ are connected
by an edge $f$ with label $x_k$. We have inequalities $0\le s<r'-i$,
$0\le t<r''-i$, where $r'$, $r''$ are the ranks of $v'$, $v''$,
respectively. The graph $Y$ is connected so there exists a path from $v'$ to
$v''$ in $Y$. The label of this path is a word $W$ that involves only letters
of the form $x_j^{\pm1}$, where $j>i$. Equalities $v'W=v''$ and
$v'x_i^{-s}x_k=v''x_i^{-t}$ imply $x_i^sW=x_kx_i^t$. We consider several
cases.

1) $k<i$. In this case $x_k$ belongs to the subgroup of $F$ generated
by $x_{k+1}$, $x_{k+2}$, \dots, which is impossible.

2) $k=i$. In this case $x_i^{t-s+1}=W$ belongs to the subgroup
generated by $x_{i+1}$, $x_{i+2}$, \dots\,. This can happen only if the
exponent $t-s+1$ equals $0$ and $W=1$ in $F$. So $v'=v''=v$ and our
vertices have the form $vx_i^{-s}$, $vx_i^{-(s+1)}$. These vertices
are connected by the edge with label $x_i$ in $\aaa_i(Y)$.

3) $k>i$. The word $W^{-1}x_i^{-s}x_kx_i^t$ equals $1$ in $F$ so the
algebraic exponent sum over all letters with the smallest subscript
must be zero. In our case, the smallest subscript is $i$ so $s=t$.
Therefore, $W=x_i^{-s}x_kx_i^s=x_{k+s}$ in $F$. The graph $Y$ is full
so the edge $e$ labelled by $x_{k+s}$ that connects vertices $v'$, $v''$
must belong to $Y$. Then our edge $f$ labelled by $x_k$ is parallel to
$e$ and thus belongs to $\aaa_i(Y)$ by definition.

Note that the rank of each vertex in $\aaa_i(Y)$ is at least $i$. This
means that we can apply the operator $\aaa_i'$ to $\aaa_i(Y)$ for
any $0\le i'<i$. We start with the ``linear" graph $\Xi_{n,m}$ that
consists of $m+1$ vertices $1$, $x_n^{-1}$, \dots, $x_n^{-m}$ and $m$
positive edges labelled by $x_n$ that connect these vertices. It is
obvious that $\Xi_{n,m}$ is full and all vertex ranks are equal to $n$.
We can apply the operator $\aaa_{n-2}$ to it ($\aaa_{n-1}$ can be also
applied but this is useless). After that, we can apply $\aaa_{n-3}$,
\dots, $\aaa_1$, $\aaa_0$. As a result we get the family of subgraphs
\be{gammas}
\Gamma_{n,m}=\aaa_0\aaa_1\cdots\aaa_{n-2}\Xi_{n,m}.
\ee
Each of them is a subgraph in $\ccc_\infty$. If we erase all edges of
$\Gamma_{n,m}$ that have labels of the form $x_k^{\pm1}$, where $k>1$,
then we get the subgraph $\bar\Gamma_{n,m}$ in $\ccc_2$ from the
statement of our Theorem.

It is possible to characterize all vertices of $\Gamma_{n,m}$ as some
words with negative exponents. We only mention that all of them will
be of the form $x_n^{-s_n}x_{n-2}^{-s_{n-2}}\cdots x_1^{-s_1}x_0^{-s_0}$,
where $s_0,s_1,\dots, s_{n-2},s_n\ge0$.

Let $Z$ be a subgraph in $\ccc_\infty$. By a {\em star\/} of a vertex $v$
in $Z$ (or a $Z$-{\em star\/}) we mean the set of labels of all edges that
come out of $v$ and belong to $Z$. Suppose that $Z'$, $Z''$ are two subsets
of the set of vertices of $Z$ and $\lambda\colon Z'\to Z''$ is a bijection.
If the $Z$-star of each vertex $v\in Z'$ coincides with the $Z$-star of
its image $\lambda(v)$, then we say that $\lambda$ is a {\em local
isomorphism\/} (in $Z$). Now let $Y$ be a finite connected full subgraph
in $\ccc_\infty$ such that all its vertices have rank $>i$. For any subset
$Y'$ of $Y$ we can denote by $\aaa_i(Y')$ the union of columns of all
vertices from $Y'$. (Note that each full subgraph is determined uniquely by
the set of its vertices so this notation does not lead to a confusion. For
a single vertex $v$ we shall write $\aaa_i(v)$ instead of $\aaa_i(\{v\})$.)
It is worth noting that if two vertices $v,w\in Y$ have the same $Y$-star,
then they have the same rank $r$ and there exists a natural bijection
between their columns, $vx_i^{-k}\mapsto wx_i^{-k}$, where $0\le k<r$.
Obviously, the $\aaa_i(Y)$-stars of $vx_i^{-k}$ and $wx_i^{-k}$ are also the
same. So we can conclude that if $Y'$ and $Y''$ are disjoint subsets and
$\lambda$ is a local isomorphism between $Y'$ and $Y''$ in $Y$, then the sets
$\aaa_i(Y')$ and $\aaa_i(Y'')$ will be also disjoint and there exists a
local isomorphism between them in $\aaa_i$.

Now we can consider two vertices $x_n^{-s}$, $x_n^{-t}$ in $\Xi_{n,m}$,
where $0<s<t<m$. Obviously, the sets $\{x_n^{-s}\}$ and $\{x_n^{-t}\}$
are locally isomorphic disjoint subsets in $\Xi_{n,m}$. Applying the
argument from the above paragraph, we see that the sets $V_s$ and $V_t$
are disjoint locally isomorphic subsets in $\Gamma_{n,m}$, where
$V_k=\aaa_0\aaa_1\cdots\aaa_{n-2}(\Xi_{n,m})$ by definition for any
$0\le k\le m$. It is clear that all the sets $V_0$, $V_1$, \dots, $V_m$
form the disjoint subdivision of the set of vertices of $\Gamma_{n,m}$.
These sets have exactly the same number of vertices. Let $\rho_k$
($0\le k\le m$) be the average degree of a vertex in the subset
$V_k$, that is, the number $\sum_{v\in V_k}\deg(v)/\#V_k$ (the degree is
taken in the whole Cayley graph). We have $\rho_1=\cdots=\rho_{m-1}=\rho$.
The density of $\Gamma_{n,m}$ will be equal to
\be{rhoo}
\bdelta(\Gamma_{n,m})=\frac{\rho_0+\rho_1+\cdots+\rho_{m-1}+\rho_m}{m+1}.
\ee
Since $\rho_0$ and $\rho_m$ are bounded, the limit of~(\ref{rhoo}) as
$m$ approaches $\infty$, is exactly $\rho$. Thus to complete the proof,
we need to show that $\rho=6(n-1)/(2n-1)$.

To calculate $\rho$, it is convenient to introduce the new family $\Gamma_n$
of auxiliary finite graphs (they are no longer subgraphs of $\ccc_\infty$).
But they will be still labelled graphs (automata) and it will be clear
from the definition of them that the set of vertices of $\Gamma_n$ is in
a bijection with the set of vertices of $V_k$ for any $1\le k<m$. Moreover,
each vertex in $\Gamma_n$ have the same star as the corresponding vertex
in $V_k$. Thus $\rho$ will be equal to $\bdelta(\Gamma_n)$. The idea is
to extend the notion of $\aaa_i$. Let $Y$ be any labelled finite graph,
where all edges have labels of the form $x_j^{\pm1}$, $j\ge0$. The rank
of a vertex is defined in the same way. Suppose that all vertices in
$Y$ have rank $>i$. For each vertex $v$ of rank $r$ we consider the
set of $r-i$ vertices $v_0$, $v_1$, \dots, $v_{r-i-1}$, where $v_0=v$.
These vertices form the column of $v$. We connect them by $r-i-1$
directed edges labelled by $x_i$, where the $k$th edge ($1\le k<r-i$)
goes from $v_k$ to $v_{k-1}$. (It is easy to see that $v_k$ here is
an analog of $vx_i^{-k}$ in the above construction.) Then for any edge
$e$ of $Y$ that has label $x_j$ and goes from $v$ to $w$, we consider $j-i$
edges $e_0$, $e_1$, \dots, $e_{j-i-1}$. Namely, for any $0\le k<j-i$,
we connect $v_k$ and $w_k$ by an edge labelled by $x_{j-k}$. Clearly,
$e_0=e$. So we get a graph that contains $Y$ and we denote it by
$\aaa_i(Y)$, as above. It is easy to see that we have an extension of
the above concept. Also it is clear that we can repeat applications of
the operators $\aaa_i$ with decreasing subscripts. If we start with
the graph $\Xi_n$ that has a single vertex and a loop labelled by $x_n$
at this vertex, then we get the graph
$$
\Gamma_{n}=\aaa_0\aaa_1\cdots\aaa_{n-2}\Xi_{n}.
$$
By $\bar\Gamma_n$ we denote the graph obtain from $\Gamma_n$ by erasing
all edges labelled by $x_i^{\pm1}$, where $i\ge2$. Since the vertex of
$\Xi_n$ has the same star as the vertex $x_n^{-k}$ in $\Xi_{n,m}$ for any
$1\le k<m$, we easily conclude from our definitions that
$\delta(\bar\Gamma_n)=\rho_k=\rho$.

The graphs $\Gamma_n$ ($n\ge1$) can be easily drawn explicitly. If a
vertex $v$ has a loop at $v$ labelled by $x_m$, then it will be clear that
the rank of $v$ equals $m$. So we can draw this vertex as a circle with
the number $m$ inside. If $Y$ is a labelled graph with labels of the
form $x_j$ ($j\ge0$), then by $\Psi(Y)$ we denote the graph obtained
from $Y$ by increasing all subscripts of the labels by $1$. We know
that $\Gamma_1$ is a single loop labelled by $x_1$. To obtain $\Gamma_{n+1}$
from $\Gamma_n$ ($n\ge1$), one has to apply $\Psi$ to $\Gamma_n$ and then
apply $\aaa_0$. This is true because $\Gamma_{n+1}=\aaa_0(\aaa_1\cdots
\aaa_{n-1}\Xi_{n+1})=\aaa_0(\Psi(\aaa_0\cdots\aaa_{n-2}\Xi_{n-1}))=
\aaa_0(\Psi(\Gamma_n))$. This gives an easy way to imagine how these
graphs look like. We illustrate this process by the following picture
that shows $\Gamma_n$ for $1\le n\le4$.

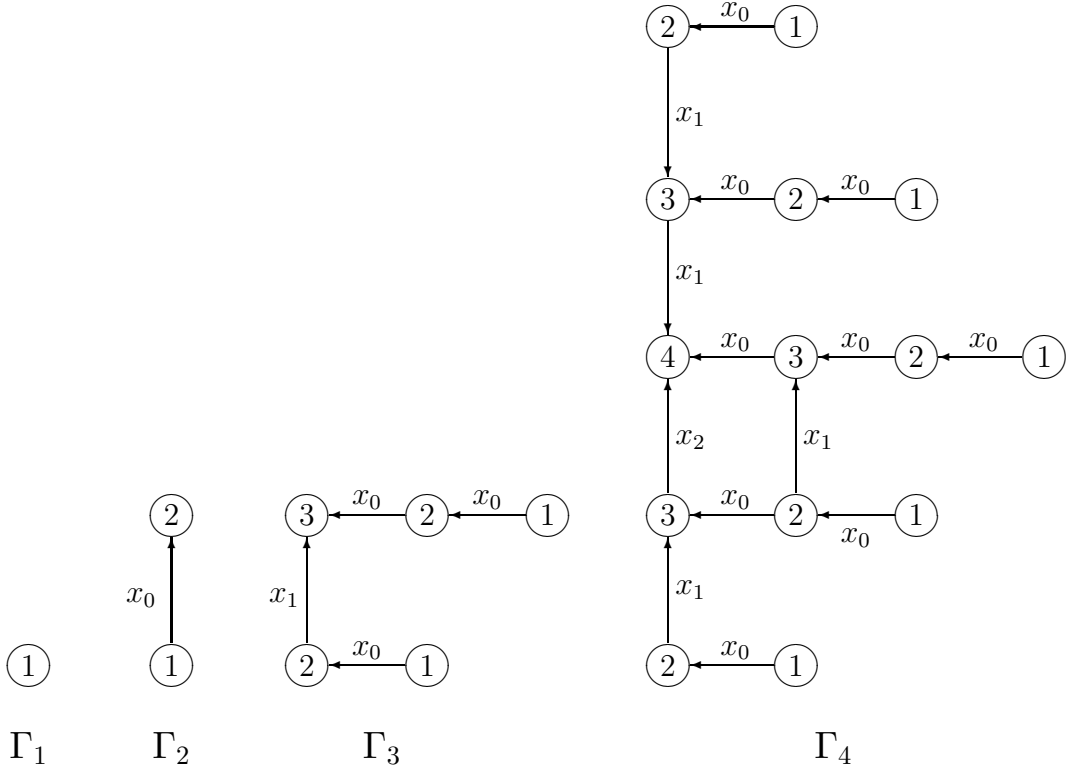
\begin{center}
\unitlength=1mm
\special{em:linewidth 0.4pt}
\linethickness{0.4pt}
\begin{picture}(141.83,100.00)
\put(4.00,12.00){\circle{6.00}}
\put(4.00,12.00){\makebox(0,0)[cc]{1}}
\put(4.00,1.00){\makebox(0,0)[cc]{\large$\Gamma_1$}}
\put(23.00,12.00){\circle{6.00}}
\put(23.00,32.00){\circle{6.00}}
\put(23.00,15.00){\vector(0,1){14.00}}
\put(23.00,12.00){\makebox(0,0)[cc]{$1$}}
\put(23.00,32.00){\makebox(0,0)[cc]{$2$}}
\put(19.00,21.00){\makebox(0,0)[cc]{$x_0$}}
\put(41.00,12.00){\circle{6.00}}
\put(41.00,32.00){\circle{6.00}}
\put(57.00,12.00){\circle{6.00}}
\put(57.00,32.00){\circle{6.00}}
\put(41.00,15.00){\vector(0,1){14.00}}
\put(54.00,12.00){\vector(-1,0){10.00}}
\put(54.00,32.00){\vector(-1,0){10.00}}
\put(41.00,12.00){\makebox(0,0)[cc]{$2$}}
\put(41.00,32.00){\makebox(0,0)[cc]{$3$}}
\put(73.00,32.00){\circle{6.00}}
\put(70.00,32.00){\vector(-1,0){10.00}}
\put(57.00,12.00){\makebox(0,0)[cc]{$1$}}
\put(57.00,32.00){\makebox(0,0)[cc]{$2$}}
\put(73.00,32.00){\makebox(0,0)[cc]{$1$}}
\put(49.00,14.00){\makebox(0,0)[cc]{$x_0$}}
\put(49.00,34.00){\makebox(0,0)[cc]{$x_0$}}
\put(65.00,34.00){\makebox(0,0)[cc]{$x_0$}}
\put(38.00,21.00){\makebox(0,0)[cc]{$x_1$}}
\put(23.00,1.00){\makebox(0,0)[cc]{\large$\Gamma_2$}}
\put(51.00,1.00){\makebox(0,0)[cc]{\large$\Gamma_3$}}
\put(89.00,12.00){\circle{6.00}}
\put(89.00,32.00){\circle{6.00}}
\put(89.00,53.00){\circle{6.00}}
\put(89.00,74.00){\circle{6.00}}
\put(89.00,97.00){\circle{6.00}}
\put(89.00,97.00){\makebox(0,0)[cc]{$2$}}
\put(89.00,74.00){\makebox(0,0)[cc]{$3$}}
\put(89.00,53.00){\makebox(0,0)[cc]{$4$}}
\put(89.00,32.00){\makebox(0,0)[cc]{$3$}}
\put(89.00,12.00){\makebox(0,0)[cc]{$2$}}
\put(89.00,15.00){\vector(0,1){14.00}}
\put(89.00,35.00){\vector(0,1){15.00}}
\put(89.00,94.00){\vector(0,-1){17.00}}
\put(89.00,71.00){\vector(0,-1){15.00}}
\put(92.00,22.00){\makebox(0,0)[cc]{$x_1$}}
\put(92.00,42.00){\makebox(0,0)[cc]{$x_2$}}
\put(92.00,64.00){\makebox(0,0)[cc]{$x_1$}}
\put(92.00,85.00){\makebox(0,0)[cc]{$x_1$}}
\put(106.00,12.00){\circle{6.00}}
\put(106.00,32.00){\circle{6.00}}
\put(106.00,53.00){\circle{6.00}}
\put(106.00,74.00){\circle{6.00}}
\put(106.00,97.00){\circle{6.00}}
\put(122.00,32.00){\circle{6.00}}
\put(122.00,53.00){\circle{6.00}}
\put(122.00,74.00){\circle{6.00}}
\put(139.00,53.00){\circle{6.00}}
\put(106.00,12.00){\makebox(0,0)[cc]{$1$}}
\put(106.00,32.00){\makebox(0,0)[cc]{$2$}}
\put(106.00,53.00){\makebox(0,0)[cc]{$3$}}
\put(106.00,74.00){\makebox(0,0)[cc]{$2$}}
\put(106.00,97.00){\makebox(0,0)[cc]{$1$}}
\put(122.00,32.00){\makebox(0,0)[cc]{$1$}}
\put(122.00,53.00){\makebox(0,0)[cc]{$2$}}
\put(122.00,74.00){\makebox(0,0)[cc]{$1$}}
\put(139.00,53.00){\makebox(0,0)[cc]{$1$}}
\put(103.00,12.00){\vector(-1,0){11.00}}
\put(103.00,32.00){\vector(-1,0){11.00}}
\put(103.00,53.00){\vector(-1,0){11.00}}
\put(103.00,74.00){\vector(-1,0){11.00}}
\put(103.00,97.00){\vector(-1,0){11.00}}
\put(119.00,32.00){\vector(-1,0){10.00}}
\put(119.00,53.00){\vector(-1,0){10.00}}
\put(119.00,74.00){\vector(-1,0){10.00}}
\put(136.00,53.00){\vector(-1,0){11.00}}
\put(106.00,35.00){\vector(0,1){15.00}}
\put(98.00,14.00){\makebox(0,0)[cc]{$x_0$}}
\put(98.00,34.00){\makebox(0,0)[cc]{$x_0$}}
\put(98.00,55.00){\makebox(0,0)[cc]{$x_0$}}
\put(98.00,76.00){\makebox(0,0)[cc]{$x_0$}}
\put(98.00,99.00){\makebox(0,0)[cc]{$x_0$}}
\put(114.00,29.00){\makebox(0,0)[cc]{$x_0$}}
\put(114.00,55.00){\makebox(0,0)[cc]{$x_0$}}
\put(114.00,76.00){\makebox(0,0)[cc]{$x_0$}}
\put(131.00,55.00){\makebox(0,0)[cc]{$x_0$}}
\put(109.00,42.00){\makebox(0,0)[cc]{$x_1$}}
\put(111.00,1.00){\makebox(0,0)[cc]{\large$\Gamma_4$}}
\end{picture}
\end{center}

\begin{lm}
\label{ank}
Let $a_{nk}$ $(1\le k\le n)$ be the number of vertices of $\Gamma_n$
that have rank $k$. Then
\be{recank}
a_{nk}=\frac{k(2n-k-1)!}{(n-k)!\,n!}.
\ee
The total number of vertices in $\Gamma_n$ equals the $n$th Catalan
number, that is,
$$
\frac{(2n)!}{n!\,(n+1)!}\,.
$$
\end{lm}

\prf We proceed by induction on $n$. Obviously, $a_{11}=1$, which agrees
with the formula. Each vertex of rank $k$ in $\Gamma_n$, where $1\le k\le n$,
becomes a vertex of rank $k+1$ in $\Psi(\Gamma_n)$. After applying $\aaa_0$
to $\Psi(\Gamma_n)$, the column of this vertex will contain exactly one
vertex of each rank from $1$ to $k+1$. Then the number of vertices of rank
$i$ in $\Gamma_{n+1}=\aaa_0(\Psi(\Gamma_n))$ will be equal to the sum
$\sum_{i\le k+1}a_{nk}$. Therefore we have the recursive formulas
$a_{n+1,1}=a_{n+1,2}=a_{n1}+\cdots+a_{nn}$, $a_{n+1,3}=a_{n2}+\cdots+a_{nn}$,
\dots, $a_{n+1,n}=a_{n,n-1}+a_{nn}$, $a_{n+1,n+1}=a_{nn}$.
By the inductive assumption, $a_{n+1,n+1}=a_{nn}=1$. Suppose that we have
already proved formula~(\ref{recank}) for $a_{n,k+1}$, where $2\le k\le n$.
Then
\begin{eqnarray*}
a_{n+1,k}&=&a_{n,k-1}+a_{n+1,k+1}=\frac{(k-1)(2n-k)!}{(n-k+1)!\,n!}+
\frac{(k+1)(2n-k)!}{(n-k)!\,(n+1)!}\\
&=&\frac{(2n-k)!}{(n-k)!\,n!}\left(\frac{k-1}{n-k+1}+\frac{k+1}{n+1}\right)=
\frac{k(2n-k+1)!}{(n-k+1)!\,(n+1)!},
\end{eqnarray*}
which proves~(\ref{recank}) for $a_{n+1,k}$. Finally,
$$
a_{n+1,1}=a_{n+1,2}=\frac{2(2n-1)!}{(n-1)!\,(n+1)!}=\frac{(2n)!}{n!\,(n+1)!},
$$
which proves~(\ref{recank}) for $a_{n+1,1}$.

Note that the total number of vertices in $\Gamma_n$ will be equal to
$a_{n1}+\cdots+a_{nn}=a_{n+1,1}=(2n)!/(n!\,(n+1)!)$, which is the $n$th
Catalan number. This completes the proof of the Lemma.
\vspace{1ex}

Now let $b_{nk}$ will be the number of edges in $\Gamma_k$ labelled
by $x_k$, where $0\le k\le n$. (Recall that if we denote a vertex by
a circle with the number $m$ inside, then this vertex has a loop
labelled by $x_m$.)

\begin{lm}
\label{bnoi}
For any $n\ge2$,
$$
b_{n0}=b_{n1}=\frac{3(2n-2)!}{(n-2)!\,(n+1)!}.
$$
\end{lm}

\prf (Remark that $b_{10}=0$, $b_{11}=1$.) We would like to establish a
more general formula
\be{bnk}
b_{nk}=\frac{(k+1)(2n-k-2)!}{(n-k)!\,(n+1)!}\cdot(3n^2-3n(k+1)+k^2+2k)
\ee
for any $0\le k<n$ together with $b_{nn}=1$. The conclusion of the
Lemma is a partial case of these. We proceed by induction on $n$. The
case $n=1$ is obvious. First of all, we want to give recursive formulas
for $b_{n+1,i}$ ($0\le i\le n+1$). Let $v$ be a vertex of rank $k$ in
$\Gamma_n$ ($1\le k\le n$). It becomes a vertex of rank $k+1$ in
$\Psi(\Gamma_n)$ and thus creates exactly $k$ edges labelled by $x_0$
in $\aaa_0(\Psi(\Gamma_n))=\Gamma_{n+1}$. This implies
$b_{n+1,0}=a_{n1}+2a_{n2}+\cdots+na_{nn}$. Using the recursive formulas
from the proof of Lemma~\ref{ank}, one can rewrite this as
$b_{n+1,0}=(a_{n1}+\cdots+a_{nn})+(a_{n2}+\cdots+a_{nn})+\cdots+a_{nn}=
a_{n+1,2}+a_{n+1,3}+\cdots+a_{n+1,n+1}=a_{n+2,3}$. Therefore,
$b_{n+1,0}=3(2n)!/((n-3)!n!)$, which proves~(\ref{bnk}) for $b_{n+1,0}$.

Let $e$ be an edge labelled by $x_k$ in $\Gamma_n$ ($0\le k\le n$). It
becomes and edge labelled by $x_{k+1}$ in $\Psi(\Gamma_n)$ so it creates one
edge of each of the ranks from $1$ to $k+1$ in $\aaa_0(\Psi(\Gamma_n))=
\Gamma_{n+1}$. Thus $b_{n+1,i}$ will be the sum $\sum_{i\le k+1}b_{nk}$
for any $1\le i\le n+1$. We have the recursive formulas
$b_{n+1,1}=b_{n0}+\cdots+b_{nn}$, $b_{n+1,2}=b_{n1}+\cdots+b_{nn}$,
\dots, $b_{n+1,n+1}=b_{nn}$. Therefore, $b_{n+1,n+1}=b_{nn}=1$ by the
inductive assumption. Suppose that~(\ref{bnk}) has been proved for
$b_{n+1,k+1}$, where $1\le k\le n$. We have $b_{n+1,k}=b_{n,k-1}+b_{n+1,k+1}$.
Hence, using the assumptions, we obtain that $b_{n+1,k}$ equals
$$
\frac{k(2n-k-1)!}{(n-k+1)!\,(n+1)!}\cdot(3n^2-3nk+k^2-1)+
\frac{(k+2)(2n-k-1)!}{(n-k)!\,(n+2)!}\cdot(3n^2-3nk+k^2+k)
$$
$$
=\frac{(2n-k-1)!}{(n-k)!\,(n+1)!}\left(\frac{k(3n^2-3nk+k^2-1)}{n-k+1}+
\frac{(k+2)(3n^2-3nk+k^2+k)}{n+2}\right)
$$
$$
=\frac{(k+1)(2n-k)!}{(n-k+1)!\,(n+2)!}\cdot(3n^2-3nk+3n+k^2-k),
$$
which proves~(\ref{bnk}) for $b_{n+1,k}$.

The proof is complete.
\vspace{1ex}

Finally, to compute $\delta(\bar\Gamma_n)$, one has to divide the number
of oriented edges in $\bar\Gamma_n$ labelled by $x_0^{\pm1}$, $x_1^{\pm1}$
by the number of vertices in $\Gamma_n$ (it is the same as the number of
vertices in $\bar\Gamma_n$). According to Lemmas~\ref{ank} and \ref{bnoi},
we have
$$
\delta(\bar\Gamma_n)=\frac{2(b_{n0}+b_{n1})}{a_{n+1,1}}=
\frac{12(2n-2)!}{(n-2)!\,(n+1)!}\cdot\frac{n!\,(n+1)!}{(2n)!}=
\frac{6(n-1)}{2n-1}
$$
for any $n\ge2$.

This completes the proof. It may be interesting to know how many
vertices in $\bar\Gamma_n$ have a given degree. For any $n\ge5$, one has
$$
\nu_2=\frac{3(2n-4)!}{(n-2)!\,(n-1)!},\
\nu_3=\frac{4(5n-12)(2n-5)!}{(n-3)!\,n!},\
\nu_4=\frac{6(2n-5)!}{(n-5)!\,(n+1)!}.
$$
Here $\nu_d$ denotes the number of vertices that have degree $d$.
This implies that an average vertex of $\bar\Gamma_{m,n}$ has degree
$2$, $3$, $4$ with probabilities approaching $3/16$, $5/8$, $3/16$,
respectively, as $m$, $n$ approach infinity.

\section{The Lower Bound for the Growth Rate}
\label{lngth}

By a {\em norm\/} of an element of $F$ we mean the shortest length of
a word in generators $\{x_0^{\pm1},x_1^{\pm1}\}$ that represents this
element. Equivalently, the norm of an element is the distance from this
element to the identity in the Cayley graph $\ccc_2$ of the group $F$.

Let $b_n$ ($n\ge0$) be the number of elements in the ball around $1$
in $\ccc_2$ of radius $n$. Since $F$ has a free subsemigroup generated
by $x_0$, $x_1$, the group $F$ has exponential growth. But it is still
unknown what will be the general growth function
$$
G(t)=\sum\limits_{n=0}^{\infty}b_n=b_0+b_1t+b_2t^2+\cdots+b_nt^n+\cdots
$$
of the group $F$ in canonical generators $x_0$, $x_1$. Since the
sequence $b_n$ is monotone and satisfy the submultiplicative inequality
$b_{m+n}\le b_mb_n$ ($m,n\ge0$), there exists the limit
$$
b=\lim\limits_{n\to\infty}\sqrt[n]{b_n},
$$
which is called the {\em growth rate\/} of a group (in a given set
of generators). The exact value of the growth rate of $F$ in generators
$x_0$, $x_1$ is also unknown. Recently Burillo~\cite{Bur} found an exact
number of {\bf positive} elements in the ball of radius $n$ (an element
of $F$ is positive if it is a product of positive letters $x_0$, $x_1$,
$x_2$, \dots of the infinite set of generators). He also found the general
growth function for positive elements and its growth rate. This implies the
lower bound for the growth rate of $F$. Namely, Burillo proved that $b$ is
not less than the largest root of the equation $x^3-x^2-2x+1=0$, that is,
$b\ge2.2469796\dots$\,. We would like to present a better lower bound
for the growth rate of $F$.

\begin{thm}
\label{sqrt5}
The growth rate of the group $F$ in generators $x_0$, $x_1$ is not
less than $(3+\sqrt5)/2=2.6180339\dots$.
\end{thm}

\prf We use the normal form for the elements of $F$ found in~\cite{GuSa}.
It was shown that there exists a regular spanning tree of the Cayley
graph $\ccc_2$. Each element in $F$ can be uniquely represented by a word
of the form $w(x_0,x_1)$, where $w$ does not contain the following forbidden
subwords: $x_i^{\pm1}x_i^{\mp1}$ ($i=0,1$), $x_1^{\pm1}x_0^nx_1$,
$x_1^{\pm1}x_0^{n+1}x_1^{-1}$ ($n\ge1$). This normal form of an element
does not give a minimal representative of an element in the sense of
its norm (that is, the regular tree is not geodesic). However, for
the regular language $L$ of normal forms one can find the growth rate using
standard methods. If a word of length $n$ belongs to this language,
then the norm of the corresponding element of $F$ does not exceed $n$
and so it will belong to $B_n(1)$. Thus $b_n$, the number of elements in
$B_n(1)$, will be not less than the number of words in $L$ that have
length $n$. So we begin to calculate the growth function of the regular
language $L$. To draw its generating automaton, we need to subdivide $L$
into seven disjoint subsets.
\vspace{1ex}

1) The empty word.

2) The words of the form $x_0^k$, where $k\ge1$.

3) The words that end with the letter $x_0^{-1}$.

4) The words that end with the letter $x_1$.

5) The words that end with the letter $x_1^{-1}$.

6) The words that end with $x_1x_0$ or $x_1^{-1}x_0$.

7) The words that end with $x_1x_0^2$ or $x_1^{-1}x_0^2$.

It is obvious that all the seven classes are pair-wise disjoint and their
union is $L$. This allows to draw the generating automaton of $L$:

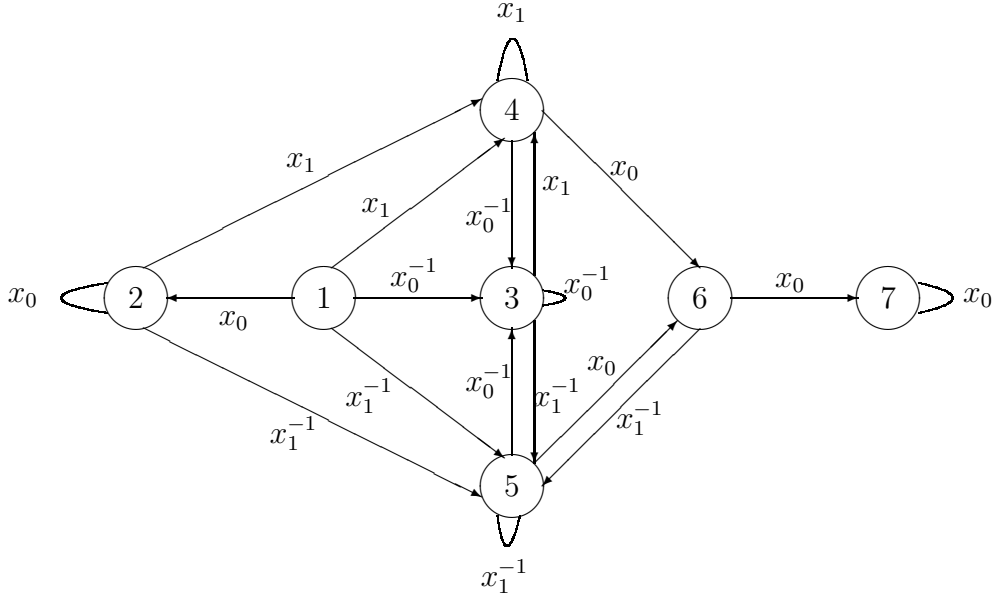
\begin{center}
\unitlength=1mm
\special{em:linewidth 0.4pt}
\linethickness{0.4pt}
\begin{picture}(130.00,78.00)
\put(67.00,38.00){\circle{8.94}}
\put(42.00,38.00){\circle{8.94}}
\put(17.00,38.00){\circle{8.94}}
\put(67.00,63.00){\circle{8.94}}
\put(67.00,13.00){\circle{8.94}}
\put(92.00,38.00){\circle{8.94}}
\put(117.00,38.00){\circle{8.94}}
\put(42.00,38.00){\makebox(0,0)[cc]{1}}
\put(17.00,38.00){\makebox(0,0)[cc]{2}}
\put(67.00,38.00){\makebox(0,0)[cc]{3}}
\put(67.00,63.00){\makebox(0,0)[cc]{4}}
\put(67.00,13.00){\makebox(0,0)[cc]{5}}
\put(92.00,38.00){\makebox(0,0)[cc]{6}}
\put(117.00,38.00){\makebox(0,0)[cc]{7}}
\put(46.00,38.00){\vector(1,0){17.00}}
\put(38.00,38.00){\vector(-1,0){17.00}}
\put(43.00,42.00){\vector(4,3){23.00}}
\put(43.00,34.00){\vector(4,-3){23.00}}
\put(18.00,34.00){\vector(2,-1){45.00}}
\put(18.00,42.00){\vector(2,1){45.00}}
\put(67.00,59.00){\vector(0,-1){17.00}}
\put(67.00,17.00){\vector(0,1){17.00}}
\put(70.00,41.00){\vector(0,1){19.00}}
\put(70.00,35.00){\vector(0,-1){19.00}}
\put(71.00,63.00){\vector(1,-1){21.00}}
\put(92.00,34.00){\vector(-1,-1){21.00}}
\put(96.00,38.00){\vector(1,0){17.00}}
\put(70.00,16.00){\vector(1,1){19.00}}
\bezier{88}(65.00,67.00)(67.00,78.00)(69.00,67.00)
\bezier{64}(65.00,9.00)(66.00,1.00)(68.00,9.00)
\bezier{72}(121.00,40.00)(130.00,38.00)(121.00,36.00)
\bezier{96}(13.00,40.00)(1.00,38.00)(13.00,36.00)
\bezier{48}(71.00,39.00)(77.00,38.00)(71.00,37.00)
\put(39.00,56.00){\makebox(0,0)[cc]{$x_1$}}
\put(38.00,20.00){\makebox(0,0)[cc]{$x_1^{-1}$}}
\put(82.00,55.00){\makebox(0,0)[cc]{$x_0$}}
\put(84.00,22.00){\makebox(0,0)[cc]{$x_1^{-1}$}}
\put(79.00,29.00){\makebox(0,0)[cc]{$x_0$}}
\put(104.00,40.00){\makebox(0,0)[cc]{$x_0$}}
\put(129.00,38.00){\makebox(0,0)[cc]{$x_0$}}
\put(67.00,76.00){\makebox(0,0)[cc]{$x_1$}}
\put(66.00,1.00){\makebox(0,0)[cc]{$x_1^{-1}$}}
\put(2.00,38.00){\makebox(0,0)[cc]{$x_0$}}
\put(77.00,40.00){\makebox(0,0)[cc]{$x_0^{-1}$}}
\put(73.00,53.00){\makebox(0,0)[cc]{$x_1$}}
\put(73.00,25.00){\makebox(0,0)[cc]{$x_1^{-1}$}}
\put(64.00,27.00){\makebox(0,0)[cc]{$x_0^{-1}$}}
\put(64.00,49.00){\makebox(0,0)[cc]{$x_0^{-1}$}}
\put(49.00,50.00){\makebox(0,0)[cc]{$x_1$}}
\put(30.00,35.00){\makebox(0,0)[cc]{$x_0$}}
\put(48.00,25.00){\makebox(0,0)[cc]{$x_1^{-1}$}}
\put(54.00,41.00){\makebox(0,0)[cc]{$x_0^{-1}$}}
\end{picture}
\end{center}

Now let $A$ be the transition matrix of the automaton, that is,
$a_{ij}=1$ if there is a directed edge from the $i$th to the $j$th
vertex of the automaton and $a_{ij}=0$ otherwise. Obviously,

$$
A=\left(
\begin{array}{ccccccc}
0&1&1&1&1&0&0\\
0&1&0&1&1&0&0\\
0&0&1&1&1&0&0\\
0&0&1&1&0&1&0\\
0&0&1&0&1&1&0\\
0&0&0&0&1&0&1\\
0&0&0&0&0&0&1
\end{array}
\right)
$$

Let $c_n^{(p)}$ be the number of directed paths of length $n$ in the
automaton from the $1$st vertex to the vertex number $p$ ($n\ge0$,
$1\le p\le7$). Clearly, $c_0^{(1)}=1$, $c_0^{(p)}=0$ if $2\le p\le7$.
For any $n\ge1$, $1\le p\le7$ one has
\be{genfs}
c_n^{(p)}=c_{n-1}^{(1)}a_{1p}+c_{n-1}^{(2)}a_{2p}+\cdots+c_{n-1}^{(7)}a_{7p}.
\ee
Let $C_p(t)$ be the generating function of $c_n^{(p)}$, that is,
$$
C_p(t)=\sum\limits_{n=0}^\infty c_n^{(p)}t^n=
c_0^{(p)}+c_1^{(p)}t+\cdots+c_n^{(p)}t^n+\cdots
$$
for each $1\le p\le7$. Obviously, from~(\ref{genfs}) we obtain $C_1=1$,
$C_p(t)=t\sum_{i=1}^7C_i(t)a_{ip}$. In the matrix form, this means that
$(C_1(t),\dots,C_7(t))(E-tA)=(1,0,\dots,0)$. Thus $(C_1(t),\dots,C_7(t))$
is $(1,0,\dots,0)(E-tA)^{-1}$, that is, the first row of $(E-tA)^{-1}$.
We are interested in the function $C(t)=C_1(t)+\cdots+C_7(t)$. Finding
the inverse matrix to $E-tA$ and adding the elements in its first row,
gives us
$$
C(t)=\frac{t^3+1}{(t-1)(t^2-3t+1)}
$$
which is presented by series as
\be{grfgusa}
C(t)=1+4t+12t^2+34t^3+92t^4+244t^5+642t^6+1684t^7+4412t^8+11554t^9+\cdots.
\ee

The coefficient $c_n$ on $t^n$ shows the number of directed paths of length
$n$ in the automaton that start in the first vertex. Hence $C(t)$ is the
generating function for the regular language $L$. The radius of convergence
of the series~(\ref{grfgusa}) will be the least absolute value of a root of
the denominator, that is, it will be $(3-\sqrt5)/2$. So the growth rate of
$L$ will be the reciprocal, that is, $(3+\sqrt5)/2$. We have already
mentioned that $\sum_{i=0}^nc_n\le b_n$ for all $n\ge0$. So we have the
lower bound for the growth function of $F$ and $(3+\sqrt5)/2=2.6180339\dots$
will be the lower bound for its growth rate. This completes the proof.

To compare the growth function of $L$ with the growth function of $F$, let
us show the first terms of the series for the function $G_s(t)$, where
the coefficient on $t^n$ below shows the number $s_n=b_n-b_{n-1}$, that is,
the number of elements in the sphere of radius $n$ in the Cayley graph
of $F$. It is easy to see that the growth rates of $b_n$ and $s_n$ are
the same. The coefficients of the series below were found on a computer.
$$
G_s(t)=1+4t+12t^2+36t^3+108t^4+314t^5+906t^6+2576t^7+7280t^8+20352t^9+\cdots.
$$
One can mention that the sequence $s_n/s_{n-1}$ is decreasing for
$1\le n\le9$. Since $s_9/s_8=2.7956043\dots$, one can expect that the
growth rate of $F$ has this number as an upper bound. (Of course, this
is only a conjecture.)

\section{Representing Elements of $F$ by Diagrams}
\label{backgr}

This background Section is devoted to the description of the representation
of elements of $F$ by semigroup diagrams. We need this for the result of
the next Section. The contents of the present Section is essentially
known. Detailed information about diagram groups can be found in~\cite{GbS}.
However, we need to describe a modified version of this idea based on the
representation of $F$ by non-spherical diagrams.

First of all, let us recall the concept of a semigroup diagram and
introduce some notation. To do this, we consider the following example.
Let $\pp=\la a,b\mid aba=b,bab=a\ra$ be the semigroup presentation. It is
easy to see by the following algebraic calculation
$$
a^5=a(bab)a(bab)a=(aba)(bab)(aba)=bab=a
$$
that the words $a^5$ and $a$ are equal modulo $\pp$. The same can be
seen from the following picture

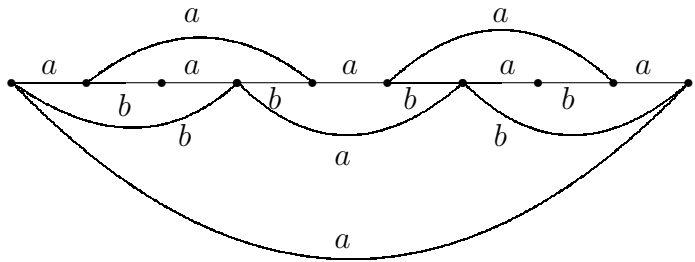
\begin{center}
\unitlength=1.00mm
\special{em:linewidth 0.4pt}
\linethickness{0.4pt}
\begin{picture}(90.00,37.00)
\put(00.00,23.00){\circle*{1.00}}
\put(10.00,23.00){\circle*{1.00}}
\put(20.00,23.00){\circle*{1.00}}
\put(30.00,23.00){\circle*{1.00}}
\put(30.00,23.00){\circle*{1.00}}
\put(40.00,23.00){\circle*{1.00}}
\put(50.00,23.00){\circle*{1.00}}
\put(60.00,23.00){\circle*{1.00}}
\put(60.00,23.00){\circle*{1.00}}
\put(70.00,23.00){\circle*{1.00}}
\put(80.00,23.00){\circle*{1.00}}
\put(90.00,23.00){\circle*{1.00}}
\emline{0.00}{23.00}{351}{90.00}{23.00}{352}
\bezier{152}(10.00,23.00)(25.00,35.00)(40.00,23.00)
\bezier{240}(50.00,23.00)(80.00,23.00)(50.00,23.00)
\bezier{164}(50.00,23.00)(65.00,37.00)(80.00,23.00)
\bezier{240}(0.00,23.00)(30.00,23.00)(0.00,23.00)
\bezier{156}(0.00,23.00)(17.00,11.00)(30.00,23.00)
\bezier{164}(30.00,23.00)(44.00,9.00)(60.00,23.00)
\bezier{164}(60.00,23.00)(74.00,9.00)(90.00,23.00)
\put(5.00,25.00){\makebox(0,0)[cc]{$a$}}
\put(24.00,32.00){\makebox(0,0)[cc]{$a$}}
\put(45.00,25.00){\makebox(0,0)[cc]{$a$}}
\put(65.00,32.00){\makebox(0,0)[cc]{$a$}}
\put(84.00,25.00){\makebox(0,0)[cc]{$a$}}
\put(23.00,16.00){\makebox(0,0)[cc]{$b$}}
\put(44.00,13.00){\makebox(0,0)[cc]{$a$}}
\put(65.00,16.00){\makebox(0,0)[cc]{$b$}}
\put(15.00,20.00){\makebox(0,0)[cc]{$b$}}
\put(24.00,25.00){\makebox(0,0)[cc]{$a$}}
\put(35.00,21.00){\makebox(0,0)[cc]{$b$}}
\put(53.00,21.00){\makebox(0,0)[cc]{$b$}}
\put(66.00,25.00){\makebox(0,0)[cc]{$a$}}
\put(74.00,21.00){\makebox(0,0)[cc]{$b$}}
\bezier{520}(0.00,23.00)(45.00,-24.00)(90.00,23.00)
\put(44.00,2.00){\makebox(0,0)[cc]{$a$}}
\end{picture}
\end{center}

This is a {\em diagram\/} $\Delta$ over the semigroup presentation $\pp$.
It is a plane graph with $10$ vertices, $15$ (geometric) edges and $6$ faces
or {\em cells\/}. Each cell corresponds to an elementary transformation of a
word, that is, a transformation of the form
$p\cdot u\cdot q\to p\cdot v\cdot q$, where $p$, $q$ are words (possibly,
empty), $u=v$ or $v=u$ belongs to the set of defining relations. The
diagram $\Delta$ has the leftmost vertex denoted by $\iota(\Delta)$ and
the rightmost vertex denoted by $\tau(\Delta)$. It also has the {\em top
path\/} $\topp(\Delta)$ and the {\em bottom path\/} $\bott(\Delta)$ from
$\iota(\Delta)$ to $\tau(\Delta)$. Each cell $\pi$ of a diagram can be
regarded as a diagram itself. The above functions $\iota$, $\tau$, $\topp$,
$\bott$ can be applied to $\pi$ as well. We do not distinguish isotopic
diagrams.

We say that $\Delta$ is a $(w_1,w_2)$-diagram whenever the label of its
top path is $w_1$ and the label of its bottom path is $w_2$. In our example,
we deal with an $(a^5,a)$-diagram. If we have two diagrams such that the
bottom path of the first of them has the same label as the top path of the
second, then we can naturally {\em concatenate\/} these diagrams by
identifying the bottom path of the first diagram with the top path of the
second diagram. The result of the concatenation of a $(w_1,w_2)$-diagram and
a $(w_2,w_3)$-diagram obviously is a $(w_1,w_3)$-diagram. We use the sign
$\circ$ for the operation of concatenation. For any diagram $\Delta$ over
$\pp$ one can consider its {\em mirror image\/} $\Delta^{-1}$. A diagram may
have {\em dipoles\/}, that is, subdiagrams of the form $\pi\circ\pi^{-1}$,
where $\pi$ is a single cell. To {\em cancel\/} (or {\em reduce\/}) the
dipole means to remove the common boundary of $\pi$ and $\pi^{-1}$ and then
to identify $\topp(\pi)$ with $\bott(\pi^{-1})$. In any digaram, we can cancel
all its dipoles, step by step. The result does not depend on the order of
cancellations. A diagram is {\em irreducible\/} whenever it has no dipoles.
The operation of cancelling dipoles has an inverse operation called the
{\em insertion\/} of a dipole. These operations induce an equivalence
relation on the set of diagrams (two diagrams are {\em equivalent\/} whenever
one can go from one of them to the other by a finite sequence of
cancelling/inserting dipoles). Each equivalence class contains exactly one
irreducible diagram.

For any nonempty word $w$, the set of all $(w,w)$-diagrams forms a monoid
with the identity element $\ve(w)$ (the diagram with no cells). The
operation $\circ$ naturally induces some operation on the set of equivalence
classes of diagrams. This operation is called a {\em product\/} and
equivalent diagram are called {\em equal\/}. (The sign $\equiv$ will be
used to denote that two diagrams are isotopic.) So the set of all
equivalence classes of $(w,w)$-diagram forms a group that is called the
{\em diagram group\/} over $\pp$ with {\em base\/} $w$. We denote this
group by ${\cal D}(\pp,w)$. We can think about this group as of the set
of all irreducible $(w,w)$-diagrams. The group operation is the
concatenation with cancelling all dipoles in the result. An inverse
element of a diagram is its mirror image. We also need one more natural
operation on the set of diagrams. By the {\em sum\/} of two diagrams we
mean the diagram obtained by identifying the rightmost vertex of the
first summand with the leftmost vertex of the second summand. This
operation is also associative. The sum of diagrams $\Delta_1$,
$\Delta_2$ is denoted by $\Delta_1+\Delta_2$.

It is known that the group $F$ is the diagram group over the simplest
semigroup presentation $\pp=\la x\mid x^2=x\ra$ with base $x$ (note that
for any base $x^k$, where $k\ge1$, we get the same group). It is not hard
to compare this representation of $F$ with other known representations. Say,
in~\cite{CFP} and many other papers, the elements of $F$ are represented by
certain pairs of rooted binary trees. If an element $g\in F$ is represented
by an (irreducible) $(x,x)$-diagram $\Delta$, then the corresponding pair of
trees can be recovered as follows. Let us consider the dual graph $\Gamma$
of the graph $\Delta$. The vertices of $\Gamma$ are the midpoints of all
edges of $\Delta$. For each $(x,x^2)$-cell of $\Delta$ we connect the
midpoint of the top edge of the cell (labelled by $x$) with the midpoints of
the edges $e'$, $e''$, where $e'e''$ is the bottom path of the cell. This
gives a ``wedge" usually called the {\em caret\/}. All these carets form a
rooted binary tree. Similarly, we do the same thing for all $(x^2,x)$-cells
of $\Delta$. This gives the other rooted binary tree. Each of the two
trees have the same number of {\em leaves\/}, that is, the vertices which
are not roots of any carets. The $i$th leaf of the upper tree (if to
count the leaves from the left to the right) coincides with the $i$th leaf
of the lower tree.

Now let us compare the diagram representation of $F$ with the
representation of its elements by piecewise-linear homeomorphisms of
the closed unit interval $[0,1]$. Let $\Delta$ be an $(x^p,x^q)$-diagram
over $\pp$. We will show how to assign to it a piecewise-linear function
from $[0,p]$ onto $[0,q]$. Each positive edge of $\Delta$ is homeomorphic
to the unit interval $[0,1]$. So we assign a coordinate to each point of
this edge (the leftmost end of an edge has coordinate $0$, the rightmost one
has coordinate $1$). Let $\pi$ be an $(x,x^2)$-cell of $\Delta$. Let us map
$\topp(\pi)$ onto $\bott(\pi)$ linearly, that it, the point on the edge
$\topp(\pi)$ with coordinate $t\in[0,1]$ is taken to the point on
$\bott(\pi)$ with coordinate $2t$ (the bottom path of $\pi$ has length
$2$ so it is naturally homeomorphic to $[0,2]$). The same thing can be
done for an $(x^2,x)$-cell of $\Delta$. Thus for any cell $\pi$ of $\Delta$
we have a natural mapping $T_\pi$ from $\topp(\pi)$ onto $\bott(\pi)$ (we
call it a {\em transition map\/}). Now let $t$ be any number in $[0,p]$. We
consider the point $o$ on $\topp(\Delta)$ that has coordinate $t$. If $o$ is
not a point of $\bott(\Delta)$, then it is an internal point on the top path
of some cell. Thus we can apply the corresponding transition map to $o$. We
repeat this operation until we get a point $o'$ on the path $\bott(\Delta)$.
The coordinate of this point is a number in $[0,q]$. Hence we have a function
$f_\Delta\colon[0,p]\to[0,q]$ induced by $\Delta$. It is easy to see this
will be a piecewise-linear function. When we concatenate diagrams, this
corresponds to the composition of the PL functions induced by these diagrams.
For the group $F$, which is the diagram group $\dd(\pp,x)$, we have the
homeomorphism from it to $PLF[0,1]$. It is known this is an isomorphism.

Now let us recall that $F$ has another representation by PL functions.
Namely, let $f_i\colon[0,\infty)\to[0,\infty)$ ($i\ge0$ is an integer)
be the PL function that has slope $1$ on $[0,i]$, slope $2$ on $[i,i+1]$
and slope $1$ again on $[i+1,\infty)$. It is easy to verify that for
any integers $j>i\ge0$, one has $f_jf_i=f_if_{j+1}$ (the functions act
on the right). Thus the mapping $x_i\mapsto f_i$ induces a homomorphism
from $F$ to $PLF[0,\infty)$. It is also known this is an isomorphism.

In the above representation of $F$ by diagrams, all these diagrams were
{\em spherical\/}, that is, they were $(w,w)$-diagrams for some word $w$.
Our new representation of $F$ by non-spherical diagrams, which we are
going to describe, corresponds to the representation of $F$ by the PL
functions on $[0,\infty)$. Let $\Delta$ be any diagram over
$\pp=\la x\mid x^2=x\ra$, not necessarily spherical. Let us add an
infinite sequence of edges on the right of $\Delta$, each edge is labelled
by $x$. This object will be called an {\em infinite diagram\/} over $\pp$.
Note that it has finitely many cells. An infinite diagram that corresponds
to $\Delta$ will be denoted by $\hat\Delta$. It has the leftmost vertex
$\iota(\hat\Delta)$ and two distinguished infinite paths starting at
$\iota(\hat\Delta)$, both labelled by the infinite power of $x$. These paths
will be denoted by $\topp(\hat\Delta)$ and $\bott(\hat\Delta)$, respectively.
The concept of a dipole in an infinite diagram is defined as above. The
same concerns the operations of deleting/inserting a dipole, the
equivalence relation induced by these operations, and so on. Any two
infinite diagrams can be naturally concatenated (the bottom path of the
first factor is identified with the top path of the second factor).
We use the sign $\circ$ for the concatenation. The operation $\circ$
gives the set of all infinite diagrams a monoid structure. The identity
of it is the infinite diagram without cells denoted by $\ve$. As in the case
of spherical diagrams, the operation of concatenation induces a group
operation on the set of all equivalence classes of infinite diagrams.
Thus we have a group. We shall denote it by $\hat\dd(\pp,x)$. (This makes
sense for any semigroup presentation $\pp$ in an alphabet of one letter.
Note that we can forget about the labels if we work with a one-letter
alphabet.) It is quite easy to see that the group we have will be isomorphic
to $F$. Indeed, let $X_i$ be the infinite diagram

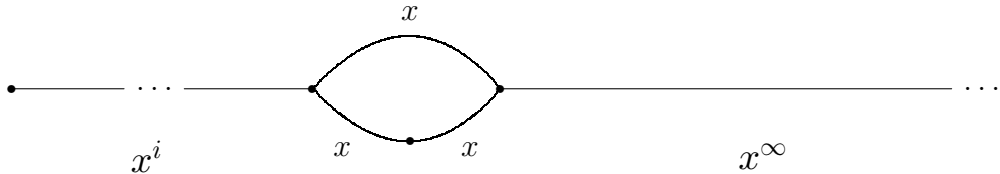
\begin{center}
\unitlength=1mm
\special{em:linewidth 0.4pt}
\linethickness{0.4pt}
\begin{picture}(131.00,25.00)
\put(2.00,11.00){\circle*{1.00}}
\emline{2.00}{11.00}{3}{17.00}{11.00}{4}
\put(21.00,11.00){\makebox(0,0)[cc]{$\dots$}}
\emline{25.00}{11.00}{9}{42.00}{11.00}{10}
\put(42.00,11.00){\circle*{1.00}}
\put(67.00,11.00){\circle*{1.00}}
\bezier{152}(42.00,11.00)(55.00,25.00)(67.00,11.00)
\bezier{152}(42.00,11.00)(55.00,-3.00)(67.00,11.00)
\emline{67.00}{11.00}{11}{127.00}{11.00}{12}
\put(55.00,4.00){\circle*{1.00}}
\put(131.00,11.00){\makebox(0,0)[cc]{$\dots$}}
\put(55.00,21.00){\makebox(0,0)[cc]{$x$}}
\put(46.00,3.00){\makebox(0,0)[cc]{$x$}}
\put(63.00,3.00){\makebox(0,0)[cc]{$x$}}
\put(20.00,2.00){\makebox(0,0)[cc]{\large$x^i$}}
\put(102.00,2.00){\makebox(0,0)[cc]{\large$x^\infty$}}
\end{picture}
\end{center}

\noindent
that consists of an $(x,x^2)$-cell, the finite path labelled by $x^i$
on the left of it and the infinite path labelled by the infinite power
of $x$, on the right. By $X_i^{-1}$ we mean the mirror image of $X_i$
under the horizontal axis symmetry. Infinite diagrams of the form
$X_i^{\pm1}$ ($i\ge0$) are called {\em atomic\/}. For any integers $j>i\ge0$,
the diagram

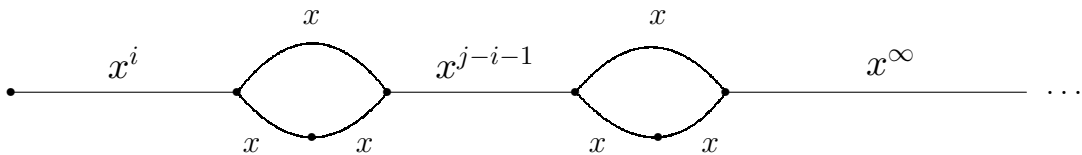
\begin{center}
\unitlength=1mm
\special{em:linewidth 0.4pt}
\linethickness{0.4pt}
\begin{picture}(142.00,22.00)
\put(2.00,9.00){\circle*{1.00}}
\put(32.00,9.00){\circle*{1.00}}
\put(52.00,9.00){\circle*{1.00}}
\put(77.00,9.00){\circle*{1.00}}
\put(97.00,9.00){\circle*{1.00}}
\emline{2.00}{9.00}{771}{32.00}{9.00}{772}
\emline{52.00}{9.00}{773}{77.00}{9.00}{774}
\emline{97.00}{9.00}{775}{137.00}{9.00}{776}
\bezier{132}(32.00,9.00)(42.00,22.00)(52.00,9.00)
\bezier{132}(32.00,9.00)(42.00,-3.00)(52.00,9.00)
\bezier{124}(77.00,9.00)(87.00,21.00)(97.00,9.00)
\bezier{124}(77.00,9.00)(88.00,-3.00)(97.00,9.00)
\put(42.00,3.00){\circle*{1.00}}
\put(88.00,3.00){\circle*{1.00}}
\put(42.00,19.00){\makebox(0,0)[cc]{$x$}}
\put(88.00,19.00){\makebox(0,0)[cc]{$x$}}
\put(34.00,2.00){\makebox(0,0)[cc]{$x$}}
\put(49.00,2.00){\makebox(0,0)[cc]{$x$}}
\put(80.00,2.00){\makebox(0,0)[cc]{$x$}}
\put(95.00,2.00){\makebox(0,0)[cc]{$x$}}
\put(17.00,13.00){\makebox(0,0)[cc]{\large$x^i$}}
\put(65.00,13.00){\makebox(0,0)[cc]{\large$x^{j-i-1}$}}
\put(119.00,13.00){\makebox(0,0)[cc]{\large$x^\infty$}}
\put(142.00,9.00){\makebox(0,0)[cc]{$\dots$}}
\end{picture}
\end{center}

\noindent
equals both $X_j\circ X_i$ and $X_i\circ X_{j+1}$. This means that we
have a homomorphism from $F$ to the group of infinite diagrams. This
homomorphism is onto because any infinite diagram is a concatenation
of atomic diagrams. The homomorphism must be injective because all
homomorphic images of $F$ are abelian~\cite{CFP}. However, the group
of infinite diagrams is not abelian since $X_1X_0=X_0X_2\ne X_0X_1$.

Note that the diagram $X_i$ corresponds to the above defined PL function
$f_i$. (The transition maps allow to assign a PL function to any infinite
diagram.) If we consider a dual graph to an infinite diagram, this
gives us a pair of binary rooted forests. Each forest is a sequence
of binary rooted trees. All of these trees but the finite number of
them consists of just the root vertex. All leaves of a forest can be
naturally enumerated. Given a pair of such forests, we can identify
the $i$th leaf of the ``upper" forest with the $i$th leaf of the ``lower"
forest for all $i\ge0$.

When we work with infinite diagrams, it is convenient to eliminate
the infinite ``tail" on the right of each infinite diagram. An ordinary
diagram over $\pp$ is called {\em canonical\/} whenever it has no dipoles
and it is not a sum of a diagram and an edge. It is obvious that there is
a one-to-one correspondence between the set of infinite diagrams without
dipoles and the set of canonical diagrams. So we may assume that each
element of $F$ has a unique canonical representative. We have a group
structure on the set of all canonical diagrams over $\pp$. Given an
$(x^p,x^q)$-diagram $\Delta_1$ and an $(x^s,x^t)$-diagram $\Delta_2$,
we multiply them as follows. If $q=s$, then we concatenate them. If
$q<s$, then we concatenate $\Delta_1+\ve(x^{s-q})$ and $\Delta_2$. If
$q>s$, then we concatenate $\Delta_1$ and $\Delta_2+\ve(x^{q-s})$.
After the concatenation, we reduce all dipoles in the result. Then we
need to make the diagram canonical. This means that we have to delete
the rightmost common suffix of the top and the bottom path of the
diagram we have. The only exception is made for the diagram $\ve(x)$,
the identity element of the group. This diagram is already canonical
so we leave it as it is.

This representation of $F$ by canonical diagrams is very convenient for
many reasons. Say, the number of cells in the canonical diagram of an
element always equals the length of the normal form of this element. Given
a normal form, it is very easy to draw the corresponding diagram, and vice
versa. The following example illustrates the diagram that corresponds to
the element $g=x_0^2x_1x_6x_3^{-1}x_0^{-2}$ represented by its normal form:

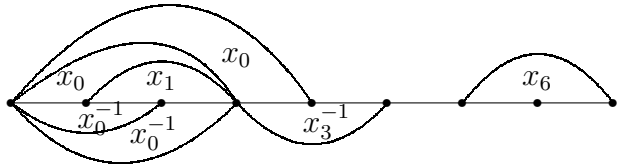
\begin{center}
\unitlength=1.00mm
\special{em:linewidth 0.4pt}
\linethickness{0.4pt}
\begin{picture}(87.00,30.00)
\put(6.00,9.00){\circle*{1.00}}
\put(16.00,9.00){\circle*{1.00}}
\put(26.00,9.00){\circle*{1.00}}
\put(36.00,9.00){\circle*{1.00}}
\put(46.00,9.00){\circle*{1.00}}
\put(56.00,9.00){\circle*{1.00}}
\put(66.00,9.00){\circle*{1.00}}
\put(76.00,9.00){\circle*{1.00}}
\put(86.00,9.00){\circle*{1.00}}
\emline{6.00}{9.00}{61}{86.00}{9.00}{62}
\bezier{132}(66.00,9.00)(76.00,22.00)(86.00,9.00)
\bezier{120}(16.00,9.00)(26.00,20.00)(36.00,9.00)
\bezier{176}(36.00,9.00)(26.00,25.00)(6.00,9.00)
\bezier{264}(46.00,9.00)(28.00,35.00)(6.00,9.00)
\bezier{120}(36.00,9.00)(46.00,-2.00)(56.00,9.00)
\bezier{104}(6.00,9.00)(16.00,1.00)(26.00,9.00)
\bezier{176}(6.00,9.00)(20.00,-7.00)(36.00,9.00)
\put(36.00,15.00){\makebox(0,0)[cc]{$x_0$}}
\put(14.00,12.00){\makebox(0,0)[cc]{$x_0$}}
\put(26.00,12.00){\makebox(0,0)[cc]{$x_1$}}
\put(76.00,12.00){\makebox(0,0)[cc]{$x_6$}}
\put(25.00,5.00){\makebox(0,0)[cc]{$x_0^{-1}$}}
\put(48.00,6.50){\makebox(0,0)[cc]{$x_3^{-1}$}}
\put(18.00,7.00){\makebox(0,0)[cc]{$x_0^{-1}$}}
\end{picture}
\end{center}

This representation of elements in $F$ by canonical diagrams has some
more advantages. First of all, the number of cells in a canonical
diagram of an element $g$ is obviously the length of the normal form
of $g$, which is also the shortest length of a word in the infinite
set of generators $\{\,x_0^{\pm1},x_1^{\pm1},x_2^{\pm1},\dots\}$ that
represents $g$. Given the canonical representative $\Delta$ of $G$,
it is easy to construct an $(x,x)$-diagram that represents the same
element. One has to add an edge on the right of $\Delta$, that is, to
consider the diagram $\Delta+\ve(x)$. Let $v$ be the rightmost vertex
of this diagram. We then connect $v$ by arcs to all the vertices in
$\topp(\Delta)$ in the upper part and connect $v$ by arcs to all the
vertices of $\Delta$ in the lower part. All edges we add have label
$x$. As a result, we give an $(x,x)$-diagram that represents $g$ in
the diagram group $\dd(\pp,x)$. This diagram has ``redundant" cells
--- they are all cells that are not in $\Delta$. These cells do not
correspond to any generators so $\Delta$ looks a more ``natural"
representative of $g$.

If we prefer to work with dual graphs, then, given a pair of binary
rooted trees from~\cite{CFP} that represent an element $g$, one has
to remove all right carets from it. The rest will be the pair of
forests that is exactly the dual graph to the canonical representative
of $g$.
\vspace{1ex}

Note that the diagram in the above picture is drawn in such a way that
all the $(x,x^2)$-cells of it are above the straight line and all the
$(x^2,x)$-cells are below the straight line. This is an easy
consequence of the fact that the normal form of an element in $F$ is
a product of the form $pq^{-1}$, where $p$, $q$ are positive words
(words in generators $x_0$, $x_1$, $x_2$, \dots with no negative
exponents). The following elementary fact was essentially used several times
in~\cite{GuSa99,Gu00} and some other papers.

\begin{lm}
\label{longpath}
Let $\pp=\la X\mid\rr\ra$ be a semigroup presentation. Suppose that all
defining relations of $\pp$ have the form $a=A$, where $a\in X$, $A$ is
a word of length at least $2$. Also assume that all letters in the
left-hand sides of the defining relations are different. Then any
irreducible diagram $\Delta$ over $\pp$ is the concatenation of the form
$\Delta_1\circ\Delta_2^{-1}$, where the top path of each cell of both
$\Delta_1$, $\Delta_2$ has length $1$. The longest positive path in $\Delta$
from $\iota(\Delta)$ to $\tau(\Delta)$ is the bottom path of both $\Delta_1$
and $\Delta_2$.
\end{lm}

Note that $\la x\mid x=x^2\ra$ obviously satisfy the conditions of the Lemma.
The same concerns the presentation $\la a,b\mid a=bab, b=aba\ra$, which
was considered in the previous Section. Let us recall the idea of the proof.
Let $p$ be the longest positive path in $\Delta$ from $\iota(\Delta)$ to
$\tau(\Delta)$. It cuts $\Delta$ into two parts. It suffices to prove that
all cells in the ``upper" part correspond to the defining relations of the
form $a=A$, where $a$ is a letter, and none of them correspond to $A=a$.
Assume the contrary. Suppose that there is a cell $\pi$ in the upper part of
$\Delta$ with the top label $A$ and the bottom label $a$. The bottom path of
$\pi$ cannot be a subpath in $p$ since $p$ is chosen the longest. So the
bottom edge of $\pi$ belongs to the top path of some cell $\pi'$. The
diagram $\Delta$ has no dipoles. All letters in the left-hand sides of the
defining relations are different. So the top path of $\pi'$ cannot have
length $1$. This means that we have found a new cell in the upper part of
$\Delta$ that also corresponds to the defining relation of the form $A=a$.
Applying the same argument to $\pi'$, we get a process that never terminates.
This is impossible since the cells that appear during the process cannot
repeat. This completes the proof.

\section{The Length Formula}
\label{lenfor}

An algorithm to find the norm of an element of $F$ in generators $x_0$,
$x_1$ was obtained by Fordham~\cite{For} in his PhD. The algorithm is
based on the representation of $F$ by pairs of binary rooted trees.
Each caret of a tree belongs to one of $7$ classes defined
in~\cite{For} (only $6$ of these classes are essential). All carets in each
of the trees are enumerated. For any number, there are $36$ cases for a
pair formed by the $i$th carets of the trees, to belong to some of the
classes. To each of these $36$ cases one can assign a weight. The sum
of these numbers will be the norm of an element. By this algorithm, the
norm can be found in a polynomial time, although the algorithm itself
is quite complicated. However, it is a powerful tool to solve various
problems about $F$. For instance, Fordham's algorithm is successfully
implemented in recent papers~\cite{Bur,CT}. In this Section, we present
another algorithm to find the norm of an element in $F$. We are based
on the representation of $F$ by non-spherical diagrams described in
Section~\ref{backgr}. We classify vertices of the canonical diagram of a
given element. There are only two cases for this. To find the norm,
one has to add the number of cells in the diagram and the number of
so-called special vertices multiplied by $2$. This algorithm is very
clear and easy in use.
\vspace{1ex}

Let $\Delta$ be a diagram over $\ppi$. From Lemma~\ref{longpath} we know
that all vertices of $\Delta$ belong to the longest positive path $p$ from
$\iota(\Delta)$ to $\tau(\Delta)$. We need to classify vertices of
$\Delta$. A non-empty positive path $q$ in $\Delta$ (in particular, a single
edge) is called a {\em brigde\/} whenever $q$ is a common subpath in
both $\topp(\Delta)$ and $\bott(\Delta)$. If $q$ is a bridge of $\Delta$,
then one has a decomposition of the form $\Delta'+q+\Delta''$, where
$\Delta$, $\Delta''$ are some diagrams (possibly $\Delta''$ is empty)
and $q$ denotes the subdiagram in $\Delta$ that consists of the bridge $q$.
We say that the bridge $q$ is {\em nontrivial\/} whenever $\Delta''$ has
cells. Note that any bridge $q$, being a path, has its initial point
$\iota(q)$.

A vertex $v$ in $\Delta$ is called {\em active\/} if $v$ is either an
initial point of a cell in $\Delta$ or $v$ is an initial point of a
nontrivial bridge of $\Delta$. Since $\Delta$ is a finite graph, we
have the distance function on vertices of $\Delta$. The distance between
two vertices is the length of the shortest path in $\Delta$ that connects
these vertices. An active vertex in $\Delta$ is called {\em special\/} if
its distance from the origin $O=\iota(\Delta)$ is strictly greater than $1$.

In the following picture

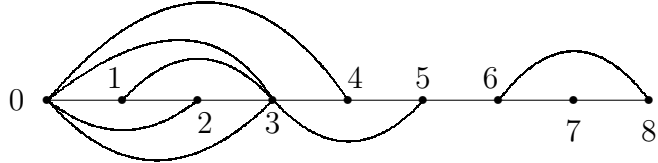
\begin{center}
\unitlength=1mm
\special{em:linewidth 0.4pt}
\linethickness{0.4pt}
\begin{picture}(87.00,30.00)
\put(6.00,9.00){\circle*{1.00}}
\put(16.00,9.00){\circle*{1.00}}
\put(26.00,9.00){\circle*{1.00}}
\put(36.00,9.00){\circle*{1.00}}
\put(46.00,9.00){\circle*{1.00}}
\put(56.00,9.00){\circle*{1.00}}
\put(66.00,9.00){\circle*{1.00}}
\put(76.00,9.00){\circle*{1.00}}
\put(86.00,9.00){\circle*{1.00}}
\emline{6.00}{9.00}{51}{86.00}{9.00}{52}
\bezier{132}(66.00,9.00)(76.00,22.00)(86.00,9.00)
\bezier{120}(16.00,9.00)(26.00,20.00)(36.00,9.00)
\bezier{176}(36.00,9.00)(26.00,25.00)(6.00,9.00)
\bezier{264}(46.00,9.00)(28.00,35.00)(6.00,9.00)
\bezier{120}(36.00,9.00)(46.00,-2.00)(56.00,9.00)
\bezier{104}(6.00,9.00)(16.00,1.00)(26.00,9.00)
\bezier{176}(6.00,9.00)(20.00,-7.00)(36.00,9.00)
\put(2.00,9.00){\makebox(0,0)[cc]{$0$}}
\put(47.00,12.00){\makebox(0,0)[cc]{$4$}}
\put(56.00,12.00){\makebox(0,0)[cc]{$5$}}
\put(65.00,12.00){\makebox(0,0)[cc]{$6$}}
\put(76.00,5.00){\makebox(0,0)[cc]{$7$}}
\put(86.00,5.00){\makebox(0,0)[cc]{$8$}}
\put(36.00,6.00){\makebox(0,0)[cc]{$3$}}
\put(27.00,6.00){\makebox(0,0)[cc]{$2$}}
\put(15.00,12.00){\makebox(0,0)[cc]{$1$}}
\end{picture}
\end{center}

\noindent
we enumerate all vertices of the diagram $\Delta$ travelling along $p$
from the origin and starting from the number $0$. It is clear that
vertices $0$, $1$, $3$, $5$, $6$ are active. Vertices $5$, $6$ from this
list are the only special vertices of $\Delta$.

\begin{thm}
\label{lngthformula}
Let $g\in F$ be represented by the canonical diagram $\Delta$. Then the
norm of $g$, that is, the length of the shortest word representing $g$
in group generators $x_0$, $x_1$, can be found by the following formula:
\be{lf}
||g||=\#_c\Delta+2\#_s\Delta,
\ee
where $\#_c$ denotes the number of cells and $\#_s$ denotes the number
of special vertices of a diagram.
\end{thm}

In the above example, one easily has $\#_c\Delta=7$, $\#_s\Delta=2$ so
$||g||=||x_0^2x_1x_6x_3^{-1}x_0^{-2}||=7+2\cdot2=11$. Recall that the number
of cells in the canonical diagram of $g$ is always equal to the length of
the normal form of $g$. Before the proof, we would like to verify that $g$
is an example of a {\em dead vertex\/} in the Cayley graph of $F$. Namely,
we would like to show that for any $a=x_i^{\pm1}$ ($i=0,1$) one has
$||ga||<||g||$ (that is, the norm of $g$ decreases in all directions). It
can be shown that the number $11$ is the smallest one for which this effect
is possible.

We need to analyze what happens if we multiply $g$ by each of the four
semigroup generators. Let us consider 4 cases.

1) $g\to gx_0$. In this transition, when we multiply $g$ by $x_0$, one
cell in $\Delta$ cancels (the leftmost cell on the bottom). The set of
vertices remains the same. It is also clear that the set of active vertices
will be the same and no vertices change their distance from the origin.
So $\#_c$ decreases by $1$, $\#_s$ is the same. Hence the norm decreases
by $1$, that is, $||gx_0||=||g||-1=10$.

2) $g\to gx_1$. When we multiply $g$ on the right by $x_1$, we cancel
the second cell on the bottom. This means we just remove the arc that
connect vertices $3$ and $5$. In this case the vertex number $3$ will be
no longer active. But the vertex number $4$ becomes active because it
will be the initial point of a bridge. However, the vertex number $4$
is on the distance $1$ from the origin so it is not special. Thus the
only special vertices in the new diagram will be $5$ and $6$, as before.
This means that the norm decreases by $1$ so $||gx_1||=||g||-1=10$.

3) $g\to gx_0^{-1}$. Multiplying by $x_0^{-1}$ on the right means
that we add a new cell to the diagram. This means we add an arc that
connects vertices $0$ and $5$. The set of active vertices remains the
same. But now the vertex number $5$ will be on the distance $1$ from
the origin so the only special vertex in the new diagram will be the
vertex number $6$. Thus $\#_c$ increases by $1$, $\#_s$ decreases by
$1$ and so the norm decreases by $1$, that is, $||gx_0^{-1}||=||g||-1=10$.

4) $g\to gx_1^{-1}$. In this case we need to add a new cell connecting
vertices $3$ and $6$ by an arc. Now the vertex number $5$ will not be
an initial point of a bridge. The only special vertex in the new diagram
will be $6$. As in the previous paragraph, $\#_c$ increases by $1$,
$\#_s$ decreases by $1$ and $||gx_1^{-1}||=||g||-1=10$.
\vspace{2ex}

{\bf Proof of Theorem~\ref{lngthformula}.}\ We use the idea
of~\cite[Lemma~{2.1.1}]{For}. For any $g\in F$ we consider the canonical
non-spherical diagram $\Delta$ that represents $g$. The value
$\#_c\Delta+2\#_s\Delta$ will be denoted by $\varphi(g)$. Our aim is to
show that $||g||=\varphi(g)$ for any $g\in F$.

Clearly, $\varphi(g)=0$ if and only if $g=1$ in $F$. We show that
$\varphi(ga)=\varphi(g)\pm1$ for any $g\in F$ and for any semigroup generator
$a\in\{x_0^{\pm1},x_1^{\pm1}\}$. This will imply that $\varphi(g)\le||g||$
for any $g\in F$. Indeed, if we decompose $g$ into a product of $n=||g||$
generators, then one can easily show by induction that $\varphi(g)\le n$.

Let us consider two neighbour elements in the Cayley graph $\ccc_2$.
For these elements we can compare the number of cells in canonical
diagrams that represent these elements. It is obvious that the difference
between the numbers of cells will be always $1$ or $-1$ so we can assume
without loss of generality that the number of cells in the canonical diagram
of $ga$ will be greater than the number of cells in the canonical diagram
$\Delta$ of $g$ by $1$.

We consider 4 cases. In each of them we need to show that the number of
special vertices either remains the same or it decreases by $1$. This
will imply that the value of $\varphi$ always increases by $1$ or decreases
by $1$ after we multiply its argument by a generator. Recall that $X_i$ is
the diagram that corresponds to the generator $x_i\in F$ ($i\ge0$).

Case 1. $a=x_0$. Since multiplying by $X_0$ increases the number of cells,
we know that $\Delta$ has no cells whose bottom path coincides with the
leftmost edge of $\bott(\Delta)$. When we concatenate $\Delta$ and $X_0$,
then a new vertex appears on the bottom path of the cell of $X_0$. This
vertex is on the distance $1$ from the origin. So the vertex we add is not
special. It is easy to see that no new special vertices can appear. So
$\varphi$ increases by $1$ in this case.

Case 2. $a=x_1$. In this case we concatenate $\Delta$ and $X_1$. This gives
a new vertex (on the bottom of the cell of $X_1$). Clearly, this vertex
is not special since it is not an initial point of a cell and it is
not the initial point of a bridge. The status of all other vertices
will be the same, that is, no special vertices appear or disappear.
Thus $\varphi$ increases by $1$.

Case 3. $a=x_0^{-1}$. Let $v_0,v_1,v_2,\dots$ be all vertices of the
bottom path of $\Delta$ enumerated from the left to the right. Multiplying
by $X_0^{-1}$ means that we connect vertices $v_0$ and $v_2$ by an
arc. Suppose that $v_2$ was a special vertex in $\Delta$. In this case
the number of special vertices in $\Delta\circ X_0^{-1}$ will be less by $1$
than the number of special vertices in $\Delta$. Otherwise the number
of special vertices remains the same. Hence $\varphi$ decreases or increases
by $1$, respectively.

Case 4. $a=x_1^{-1}$. In the notation of the previous paragraph, one has
to add an arc that connects vertices $v_1$ and $v_3$ in order to
concatenate $\Delta$ and $X_1^{-1}$. Note that vertices $v_0$, $v_1$
are not special in our diagrams. The vertex $v_3$ also does not
change its status because the new arc does not connect it with the
origin. Clearly, nothing happens with vertices $v_k$ if $k\ge4$. So
the only vertex we need to think about is $v_2$. If the distance from
it to the origin in $\Delta$ is $1$, then the same is true for
$\Delta'=\Delta\circ X_1^{-1}$. Suppose that the distance in $\Delta$ from
$v_2$ to the origin is at least $2$. Then the same is true for $\Delta'$.
If $v_2$ was not active in $\Delta$, then it could not be active in
$\Delta'$. Suppose that $v_2$ is an initial point of a cell in
$\Delta$. In this case $v_2$ be a special vertex of $\Delta$. It is
obvious that $v_2$ has the same status in $\Delta'$. So the only case
we need to check is when $v_2$ is an initial point of a nontrivial
bridge in $\Delta$. After we concatenate $\Delta$ and $X_1^{-1}$,
this will no longer be true. So $v_2$ will not be special in $\Delta'$
so the number of special vertices decreases by $1$. As a result, we
have that $\varphi$ increases or decreases by $1$.
\vspace{1ex}

To complete the proof, we need to check that for any $g\in F$, $g\ne1$,
there exists at least one generator $a\in\{x_0^{\pm1},x_1^{\pm1}\}$
such that $\varphi$ decreases at its direction, that is, $\varphi(ga)=
\varphi(g)-1$. This will imply the inequality $||g||\le\varphi(g)$ for
any $g\in F$. Indeed, if $\varphi(g)=n$, then we can start from the vertex
$g$ in $\ccc_2$. At each step we choose a direction that decreases
$\varphi$ by $1$. After $n$ steps, we get to the element with the zero
value of $\varphi$, that is, to the identity. So the distance from $g$ to
the origin (that is, the norm of $g$) does not exceed $n$.

We enumerate the vertices of $\bott(\Delta)$ as above, where $\Delta$
is the canonical representative of $g$. Since $g\ne1$, the diagram
$\Delta$ is nontrivial. Suppose that $v_2$ is a special vertex in
$\Delta$. Then $\Delta$ cannot be the concatenation of some diagram
$\Delta'$ and $X_0$ (otherwise $v_2$ is on the distance $1$ from the
origin). So if we multiply $\Delta$ by $X_0^{-1}$, then the number of
cells increases by $1$ and we are in the situation of Case 3. Since
$v_2$ is special vertex of $\Delta$, the value of $\varphi$ decreases
by $1$, that is, $\varphi(gx_0^{-1})=\varphi(g)-1$.

Assume that $v_2$ is not a special vertex of $\Delta$. Let $\Delta$
be a concatenation of some diagram $\Delta'$ and $X_1^{-1}$. If we
multiply $g$ by $x_1$ on the right, then the number of cells decreases
and we are in the situation, which is inverse to the one described in
Case 4 (now $\Delta$ is obtained from $\Delta'$ by adding a cell when
multiplying by $x_1^{-1}$). We know that $\varphi(\Delta')$ is
always greater by $1$ than $\varphi(\Delta)$ except for one case.
Let us describe it. By $\pi$ we denote the cell in $\Delta$ that
corresponds to the factor $X_1^{-1}$ in the concatenation
$\Delta'\circ X_1^{-1}$. The top path of $\pi$ is a product of two
edges $e'e''$. The terminal point of $e'$ (=the initial point of $e''$)
will be denoted by $v$. According to Case 4 (note that we have to exchange
the r\^{o}les of $\Delta$ and $\Delta'$) we have that $v$ must be the
initial point of a nontrivial bridge in $\Delta'$. This bridge starts with
the edge $e''$. So the terminal point of $e''$ must be a special vertex
of $\Delta'$ and also of $\Delta$. But this vertex is $v_2$, which
gives a contradiction.

Let us analyze all cases that can happen to the edge $e$ that is
contained in $\bott(\Delta)$ (or its continuation) and has the endpoints
$v_1$, $v_2$. Suppose that $e$ does not belong to a cell so it is a bridge.
This bridge should be trivial since $v_2$ is not special. Thus $\Delta$
is in fact an $(x^m,x)$-diagram for some $m\ge1$. It is nontrivial so
$\Delta$ is a concatenation of some diagram $\Delta'$ and $X_0^{-1}$.
If we look at the situation of Case 3 (in our case we go from $\Delta'$
to $\Delta$ adding a cell), then we see that $\varphi(\Delta)$ will
be always greater than $\varphi(\Delta')$ by $1$ except for one case
when the vertex $v_1$ is a special vertex of the subdiagram $\Delta'$.
This cannot happen because $\bott(\Delta)$ has length $1$. Thus
$\varphi(gx_0)=\varphi(g)-1$ in this case.

Now let $e$ be contained in some cell of $\Delta$. The only case we
need to consider is the one when $e$ is a part of the bottom path
of an $(x,x^2)$-cell (the case of an $(x^2,x)$-cell means that $\Delta$
is a concatenation of some diagram and $X_1^{-1}$, which has been
already considered). Let us denote this $(x,x^2)$-cell by $\pi$. The
edge $e$ may be the first or the second edge of $\bott(\pi)$ so we
have two subcases.

Suppose that $\bott(\pi)=e'e$ for some edge $e'$. If $e'$ is contained
in the bottom path of $\Delta$, then it coincides with the first edge of
$\bott(\Delta)$. Thus $\Delta$ will be the concatenation of some
diagram $\Delta'$ and $X_0$. From Case 1 we can extract that
$\varphi(gx_0^{-1})=\varphi(g)-1$ since $\Delta$ is obtained from
$\Delta'$ by adding one cell that corresponds to $X_0$. If $e'$ is not
contained in $\bott(\Delta)$, then the first edge of $\bott(\Delta)$
must be a bottom edge of an $(x^2,x)$-cell. This means that $\Delta$
is a concatenation of some diagram $\Delta'$ and $X_0^{-1}$. The
situation we have is described in Case 3. Adding the new cell to
$\Delta'$ that corresponds to $X_0^{-1}$ always increases $\varphi$
except for the case when the vertex $v_1$ will be a special vertex
of the subdiagram $\Delta'$. But we know that the only positive
edge of $\Delta$ that comes out of $v_1$ is $e$. So $v_1$ is not
an initial point of a cell. Neither $e$ is a bridge. So $\varphi(gx_0)=
\varphi(g)-1$ in this subcase.

Finally, let $\bott(\pi)=ee''$ for some edge $e''$. If $e''$ is also
contained in $\bott(\Delta)$, then $\Delta$ is a concatenation of
some diagram and $X_1$. From Case 2 we can easily conclude that
$\varphi(gx_1^{-1})=\varphi(g)-1$. If $e''$ is not contained in
$\bott(\Delta')$, then there is an ``angle" between $e''$ and the
third edge of $\bott(\Delta)$. Then $v_2$ will be an initial point
of some cell. However, $v_2$ cannot be on the distance $1$ from the
origin because otherwise the arc that connects the origin and $v_2$
would cross the edge $\topp(\pi)$, which is impossible. Therefore,
$v_2$ is special. We have a contradiction.

The proof is complete. Note that we have a quick procedure to find
minimal representatives of elements of $F$ using the above description.

\end{document}